\documentclass[12pt,reqno]{amsart}
\usepackage{amsmath,amsthm,amssymb,amscd}

\begin{document}

\newcommand{\TITLE}{$p$-adic properties of division polynomials and elliptic
divisibility sequences}
\newcommand{\TITLERUNNING}{$p$-adic properties of division polynomials}
\newcommand{\DATE}{April 2004}
\newcommand{\VERSION}{}

\theoremstyle{plain}
\newtheorem{theorem}{Theorem}
\newtheorem{conjecture}[theorem]{Conjecture}
\newtheorem{proposition}[theorem]{Proposition}
\newtheorem{lemma}[theorem]{Lemma}
\newtheorem{corollary}[theorem]{Corollary}

\theoremstyle{definition}
\newtheorem{definition}{Definition}

\theoremstyle{remark}
\newtheorem{remark}{Remark}
\newtheorem{example}{Example}
\newtheorem*{question}{Question}
\newtheorem*{acknowledgement}{Acknowledgements}

\def\BigStrut{\vphantom{$(^{(^(}_{(}$}} 

\newenvironment{notation}[0]{%
  \begin{list}%
    {}%
    {\setlength{\itemindent}{0pt}
     \setlength{\labelwidth}{4\parindent}
     \setlength{\labelsep}{\parindent}
     \setlength{\leftmargin}{5\parindent}
     \setlength{\itemsep}{0pt}
     }%
   }%
  {\end{list}}

\newenvironment{parts}[0]{%
  \begin{list}{}%
    {\setlength{\itemindent}{0pt}
     \setlength{\labelwidth}{1.5\parindent}
     \setlength{\labelsep}{.5\parindent}
     \setlength{\leftmargin}{2\parindent}
     \setlength{\itemsep}{0pt}
     }%
   }%
  {\end{list}}
\newcommand{\Part}[1]{\item[\upshape#1]}

%
\newcommand{\EndProofAtDisplay}{\renewcommand{\qedsymbol}{}}
\newcommand{\qedtag}{\tag*{\qedsymbol}}

\renewcommand{\a}{\alpha}
\renewcommand{\b}{\beta}
\newcommand{\g}{\gamma}
\renewcommand{\d}{\delta}
\newcommand{\e}{\epsilon}
\newcommand{\f}{\phi}
\newcommand{\fhat}{{\hat\phi}}
\renewcommand{\l}{\lambda}
\renewcommand{\k}{\kappa}
\newcommand{\lhat}{\hat\lambda}
\newcommand{\m}{\mu}
\renewcommand{\o}{\omega}
\renewcommand{\r}{\rho}
\newcommand{\rbar}{{\bar\rho}}
\newcommand{\s}{\sigma}
\newcommand{\sbar}{{\bar\sigma}}
\renewcommand{\t}{\tau}
\newcommand{\z}{\zeta}

\newcommand{\D}{\Delta}
\newcommand{\F}{\Phi}
\newcommand{\G}{\Gamma}

\newcommand{\ga}{{\mathfrak{a}}}
\newcommand{\gb}{{\mathfrak{b}}}
\newcommand{\gc}{{\mathfrak{c}}}
\newcommand{\gd}{{\mathfrak{d}}}
\newcommand{\gm}{{\mathfrak{m}}}
\newcommand{\gn}{{\mathfrak{n}}}
\newcommand{\gp}{{\mathfrak{p}}}
\newcommand{\gq}{{\mathfrak{q}}}
\newcommand{\gP}{{\mathfrak{P}}}
\newcommand{\gQ}{{\mathfrak{Q}}}

\def\Acal{{\mathcal A}}
\def\Bcal{{\mathcal B}}
\def\Ccal{{\mathcal C}}
\def\Dcal{{\mathcal D}}
\def\Ecal{{\mathcal E}}
\def\Fcal{{\mathcal F}}
\def\Gcal{{\mathcal G}}
\def\Hcal{{\mathcal H}}
\def\Ical{{\mathcal I}}
\def\Jcal{{\mathcal J}}
\def\Kcal{{\mathcal K}}
\def\Lcal{{\mathcal L}}
\def\Mcal{{\mathcal M}}
\def\Ncal{{\mathcal N}}
\def\Ocal{{\mathcal O}}
\def\Pcal{{\mathcal P}}
\def\Qcal{{\mathcal Q}}
\def\Rcal{{\mathcal R}}
\def\Scal{{\mathcal S}}
\def\Tcal{{\mathcal T}}
\def\Ucal{{\mathcal U}}
\def\Vcal{{\mathcal V}}
\def\Wcal{{\mathcal W}}
\def\Xcal{{\mathcal X}}
\def\Ycal{{\mathcal Y}}
\def\Zcal{{\mathcal Z}}

\renewcommand{\AA}{\mathbb{A}}
\newcommand{\BB}{\mathbb{B}}
\newcommand{\CC}{\mathbb{C}}
\newcommand{\FF}{\mathbb{F}}
\newcommand{\GG}{\mathbb{G}}
\newcommand{\NN}{\mathbb{N}}
\newcommand{\PP}{\mathbb{P}}
\newcommand{\QQ}{\mathbb{Q}}
\newcommand{\RR}{\mathbb{R}}
\newcommand{\ZZ}{\mathbb{Z}}

\def \bfa{{\mathbf a}}
\def \bfb{{\mathbf b}}
\def \bfc{{\mathbf c}}
\def \bfe{{\mathbf e}}
\def \bff{{\mathbf f}}
\def \bfF{{\mathbf F}}
\def \bfg{{\mathbf g}}
\def \bfn{{\mathbf n}}
\def \bfp{{\mathbf p}}
\def \bfr{{\mathbf r}}
\def \bfs{{\mathbf s}}
\def \bft{{\mathbf t}}
\def \bfu{{\mathbf u}}
\def \bfv{{\mathbf v}}
\def \bfw{{\mathbf w}}
\def \bfx{{\mathbf x}}
\def \bfy{{\mathbf y}}
\def \bfz{{\mathbf z}}
\def \bfX{{\mathbf X}}
\def \bfU{{\mathbf U}}
\def \bfmu{{\boldsymbol\mu}}

\newcommand{\Gbar}{{\bar G}}
\newcommand{\Kbar}{{\bar K}}
\newcommand{\kbar}{{\bar k}}
\newcommand{\Obar}{{\bar O}}
\newcommand{\Pbar}{{\bar P}}
\newcommand{\Rbar}{{\bar R}}
\newcommand{\Qbar}{{\bar Q}}
\newcommand{\QQbar}{{\bar{\QQ}}}


\newcommand{\Aut}{\operatorname{Aut}}
\newcommand{\Disc}{\operatorname{Disc}}
\renewcommand{\div}{\operatorname{div}}
\newcommand{\Div}{\operatorname{Div}}
\newcommand{\Etilde}{{\tilde E}}
\newcommand{\End}{\operatorname{End}}
\newcommand{\Fix}{\operatorname{Fix}}
\newcommand{\Frob}{\operatorname{Frob}}
\newcommand{\Gal}{\operatorname{Gal}}
\newcommand{\GCD}{{\operatorname{GCD}}}
\renewcommand{\gcd}{{\operatorname{gcd}}}
\newcommand{\hhat}{{\hat h}}
\newcommand{\Hom}{\operatorname{Hom}}
\newcommand{\Ideal}{\operatorname{Ideal}}
\newcommand{\Image}{\operatorname{Image}}
\newcommand{\longhookrightarrow}{\lhook\joinrel\relbar\joinrel\rightarrow}
\newcommand{\LS}[2]{\genfrac(){}{}{#1}{#2}}  
\newcommand{\MOD}[1]{~(\textup{mod}~#1)}
\newcommand{\Norm}{{\textup{\textsf{N}}}}
\newcommand{\NS}{\operatorname{NS}}
\newcommand{\notdivide}{\nmid}
\newcommand{\ord}{\operatorname{ord}}
\newcommand{\Parity}{\operatorname{Parity}}
\newcommand{\Pic}{\operatorname{Pic}}
\newcommand{\Proj}{\operatorname{Proj}}
\newcommand{\rank}{\operatorname{rank}}
\newcommand{\res}{\operatornamewithlimits{res}}
\newcommand{\Resultant}{\operatorname{Resultant}}
\renewcommand{\setminus}{\smallsetminus}
\newcommand{\sign}{\operatorname{Sign}}
\newcommand{\Spec}{\operatorname{Spec}}
\newcommand{\Support}{\operatorname{Support}}
\newcommand{\tors}{{\textup{tors}}}
\newcommand\W{W^{\vphantom{1}}} 
\newcommand{\Wtilde}{{\widetilde W}}
\newcommand{\<}{\langle}
\renewcommand{\>}{\rangle}

\hyphenation{para-me-tri-za-tion}


\title[\TITLERUNNING]{\TITLE}
\date{\DATE}
\author{Joseph H. Silverman}
\address{Mathematics Department, Box 1917, Brown University,
Providence, RI 02912 USA}
\email{jhs@math.brown.edu}

\thanks{The author's research supported by NSA grant H98230-04-1-0064}


\begin{abstract}
For a fixed rational point $P\in E(K)$ on an elliptic curve, we
consider the sequence of values~$\bigl(F_n(P)\bigr)_{n\ge1}$ of the
division polynomials of~$E$ at~$P$.  For a finite field~$K/\FF_p$, we
prove that the sequence is periodic.  For a local field~$K/\QQ_p$, we
prove (under certain hypotheses) that there is a power $q=p^e$ so that
for all $m\ge1$, the limit of $F_{mq^{k}}(P)$ as $k\to\infty$ exists
in~$K$ and is algebraic over~$\QQ(E)$. We apply this result to prove
an analogous $p$-adic limit and algebraicity result for elliptic
divisibility sequences.
\end{abstract}

\maketitle

\tableofcontents  

\newpage
\section*{Introduction}
Let~$E/K$ be an elliptic curve given by a Weierstrass equation
\[
  E:y^2+a_1xy+a_3y = x^3+a_2x^2+a_4x+a_6,
\]
and let $z=-x/y$ be the usual uniformizer at~$\Ocal$.
The $n$-division polynomial~$F_n$ of~$E$
is the function $F_n\in K[x,y]\subset K(E)$ with divisor
\[
  (F_n) = [n]^*(\Ocal) - n^2(\Ocal),
\]
suitably normalized at~$\Ocal$ (see
Definition~\ref{definition:ndivpoly}). Division polynomials
play  for elliptic curves the role that is played by the 
polynomials~\text{$X^n-1$} for the multiplicative group.
\par
Complex analytically the $n$-division polynomial of an elliptic
curve $\CC/L$ is equal to the quotient $\s(n\z,L)/\s(\z,L)^{n^2}$ of
Weierstrass $\s$-func\-tions~\cite[Chapter~XX, Misc{.}~examples~24
and~33]{WW}.  Division polynomials play a prominent role in the theory of
elliptic functions and elliptic curves, appearing in the elliiptic 
addition law, in the theory of complex multiplication, in 
transformation formulas for canonical local heights, in the theory
of elliptic divisibility sequences, and in the cryptographically
important SEA algorithm of Schoof, Elkies, and Atkins~\cite{Schoof1,Schoof2}
for counting points on elliptic curves over~$\FF_p$.
\par
In this paper we study the sequence of values~$(F_n(P))_{n\ge1}$ of
the division polynomials evaluated at a point~$P\in E(K)$.  We will be
especially interested in periodicity properties when~$K$ is a finite
ring or finite field and in convergence properties when~$K$ is a
complete local field. We now describe special cases of our three main
theorems for~$\FF_p$,~\text{$\ZZ/p^\mu\ZZ$}, and~$\QQ_p$. See
Corollary~\ref{corollary:divpolyperiodic},
Theorem~\ref{theorem:divpolyconverge}, and
Theorem~\ref{theorem:Fnperiodmodptothemu} for general statements.

\begin{theorem}
\label{theorem:intro:divpolyperiodic}
Let $E/\FF_p$ be an elliptic curve and let $P\in E(\FF_p)$ be a point
of exact order $r\ge2$. Then the sequence $(F_n(P))_{n\ge0}$ is
periodic with period~$rt$ for some integer~$t$ dividing~\text{$p-1$}
if~$r\ge3$ and dividing~\text{$2p-2$} if $r=2$.
\end{theorem}

Our proof of Theorem~\ref{theorem:intro:divpolyperiodic}, which is
modeled after a similar result by Ward~\cite[Theorems~8.1
and~9.2]{Ward1} for elliptic divisibility sequences, proceeds by first
lifting to a field of characteristic~$0$ and then embedding the
problem into~$\CC$ and using the transformation law for the
Weierstrass $\s$-function.

\begin{theorem}
\label{theorem:intro:divpolyconverge}
Let $E/\QQ_p$ be an elliptic curve with good ordinary reduction,
let~$P\in E(\QQ_p)$ be a point whose reduction modulo~$p$ has order~$r$,
and let~$rt$ be the period of the sequence~$(F_n(P))_{n\ge0}$ 
{\upshape(}cf. Theorem~\ref{theorem:intro:divpolyperiodic}{\upshape)}.
Assume further that~$p\ge3$, that~$r\ge3$, and that $p\notdivide r$.
Fix a power~$p^e$ of~$p$ satisfying \text{$p^e\equiv1\pmod{rt}$}. 
Then for every $m\ge1$, the limit
\[
  G_m(P) := \lim_{k\to\infty} F_{mp^{ek}}(P)\quad
  \text{converges in $\ZZ_p$.}
\]
Further, $G_m(P)=0$ if and only if~$r|m$.
\par
If in addition~$E$ is defined over~$\QQ$ and~$P\in E(\QQ)$,
then~$G_m(P)$ is algebraic over~$\QQ$.
\end{theorem}

The proof of Theorem~\ref{theorem:intro:divpolyconverge} uses the
Mazur-Tate $p$-adic $\s$-function~\cite{MT}, which is why the
statement of the theorem is restricted to the case of curves with
ordinary reduction. However, it is likely that the statement is true
in general. We will use similar techniques to prove the following
periodicity result modulo higher powers of~$p$.  (See
Remark~\ref{remark:periodicitymodppower}.) This may be compared
with Shipsey~\cite[Theorem~3.5.4]{SHIP}, who uses explicit formulas to
prove an analogous result for elliptic divisibility sequences
modulo~$p^2$.

\begin{theorem}
\label{theorem:intro:periodicitymodppower}
With notation and assumptions as in
Theorem~\ref{theorem:intro:divpolyconverge}, for every~$\mu\ge1$,
the sequence
\[
  \bigl(F_{kr}(P) \bmod p^\mu\bigr)_{k\ge1}
\]
is periodic with period dividing $p^{\mu-1}(p-1)$.
\end{theorem}

As an application of Theorem~\ref{theorem:intro:divpolyconverge}, we
will partially answer a question raised in~\cite{SS} concerning elliptic
divisibility sequences. A (proper) \emph{elliptic divisibility sequence}
is a sequence~$\Wcal=(W_n)_{n\ge0}$ of integers whose initial
terms satisfy $W_0=0$, $W_1=1$, $W_2W_3\ne0$, $W_2|W_4$, and 
whose subsequent terms are determined by the nonlinear recursion
\[
  \W_{m+n}  \W_{m-n} = \W_{m+1}  \W_{m-1}  W_n^2
     - \W_{n+1}  \W_{n-1}  W_m^2 
\]
for all $m\ge n\ge1$. Ward, who made an extensive study of these
sequences~\cite{Ward1,Ward2}, shows that a proper elliptic
divisibility sequence~$\Wcal$ is associated to a (possibly singular)
elliptic curve~$E_\Wcal$ and point~$P_\Wcal\in E_\Wcal(\QQ)$ and that
the values of~$W_n$ are closely related to the values of the division
polynomials~$F_n(P_\Wcal)$.  More recently, elliptic divisibility
sequences have been studied by Shipsey~\cite{SHIP}, who gives an
application to the elliptic curve discrete logarithm problem, and by
several other authors~\cite{CC,Durst,EGW,EPSW,EW1,EW2,SW,SS}. (See
also~\cite{Gale,Propp,Rob} for work on the related Somos sequences.)
\par
In~\cite{SS}, Nelson Stephens and the author proved that for any fixed
modulus~$2^e$, the sequence
\[
  (W_{2^k}\bmod2^e)_{k\ge0}
\]
is eventually periodic, and the question was raised as to whether this
periodicity reflects a subtler $p$-adic convergence property. 
In this paper we use the results cited above to
prove $p$-adic convergence for almost all primes.

\begin{theorem}
\label{theorem:intro:EDSconvergence}
Let $\Wcal=(W_n)_{n\ge0}$ be a proper elliptic divisibility sequence,
and assume that the associated elliptic curve~$E_\Wcal$ is nonsingular
and does not have complex multiplication.  Then for almost all primes~$p$,
in the sense of density, the following two statements are true.
\begin{parts}
\Part{(a)}
There is an exponent~$N=N_p\ge1$ so that for every $m\ge1$, the
limit
\begin{equation}
  \label{equation:padiclimit}
  \lim_{k\to\infty} W_{mp^{k\t}}\quad\text{converges in $\ZZ_p$.}
\end{equation}
\Part{(b)}
The limit given by~\eqref{equation:padiclimit} is algebraic over~$\QQ$. 
\end{parts}
{\upshape(}If~$E_\Wcal$ has complex multiplication,
then~{\upshape(}a{\upshape)} and~{\upshape(}b{\upshape)} are true for
at least half of all primes, more precisely, they are true for all but
finitely many of the primes that split in the CM~field
of~$E_\Wcal$.{\upshape)}
\end{theorem}

Although we are only able to prove
Theorem~\ref{theorem:intro:EDSconvergence} for almost all primes, we
see no reason why it should not be true in general.

\begin{conjecture}
\label{conjecture:intro:padicconvergence}
Statements~{\upshape(a)} and~{\upshape(b)} of
Theorem~\ref{theorem:intro:EDSconvergence} are true for every proper
elliptic divisibility sequence and  for every prime~$p$.
\end{conjecture}

In our proof of Theorem~\ref{theorem:intro:EDSconvergence}, it is
natural to avoid a certain finite set of primes at which~$(W_n)$
behaves badly.  However, the reason that we ultimately eliminate
infinitely many primes is because we used the Mazur-Tate $p$-adic
$\s$-function~\cite{MT} in our proof of
Theorem~\ref{theorem:intro:divpolyconverge}, so that theorem only
applies to elliptic curves with ordinary reduction.  A theorem of
Serre~\cite{SE1} says that for a fixed (non-CM) elliptic
curve~$E/\QQ$, almost all primes are ordinary, but since
Elkies~\cite{EL1} has shown that there are also infinitely many primes
of supersingular reduction~\cite{EL1}, our proof cannot be directly
extended to prove Theorem~\ref{theorem:intro:EDSconvergence} for all
but finitely many primes.

\section{Elliptic curves and division polynomials}
Let $E/K$ be an elliptic curve defined over any field~$K$,
and fix a Weierstrass equation for~$E$,
\[
  E:y^2+a_1xy+a_3y = x^3+a_2x^2+a_4x+a_6.
\]
Then~$z=-x/y\in K(E)$ is a uniformizer at~$\Ocal\in E$, and 
the invariant differential $\o=dx/(2y+a_1x+a_3)$ can be expanded
as a formal (Laurant) series in a formal neighborhood of~$\Ocal$ as
\[
  \o(z) = (1+a_1z+(a_1^2+a_2)z^2+\cdots)\,dz.
\]
This series has coefficients in~$\ZZ[a_1,a_2,a_3,a_4,a_6]$,
and we have normalized matters so that~$(\o/dz)(\Ocal)=1$.

\begin{definition}
\label{definition:ndivpoly}
Let $n\ge1$ be an integer and let $[n](z)\in K[[z]]$ be the power
series defining the multiplication-by-$n$ map on the formal group
of~$E$.  The \emph{$n$-division polynomial~$F_n$}
(normalized relative to the uniformizer~$z$) is the unique rational
function~$F_n\in K(E)$ satisfying
\[
  (F_n) = [n]^*(\Ocal) - n^2(\Ocal)
  \quad\text{and}\quad
  \left(\frac{z^{n^2}F_n}{[n](z)}\right)(\Ocal)=1.
\]
If~$n\ne 0$ in~$K$, then~$[n](z)=nz+O(z^2)$, so the normalization
condition becomes simply $(z^{n^2-1}F_n)(\Ocal)=n$. (For a more
general normalization procedure, see Remark~\ref{remark:MTdivpolydef}.)
\end{definition}

We will use the following elementary ``chain rule'' for
division polynomials.

\begin{lemma}
\label{lemma:divpolychainrule}
For all integers $m,n\ge1$, 
\[
  F_{mn} = \bigl(F_n\circ[m]\bigr)\circ F_m^{n^2}.
\]
\end{lemma}
\begin{proof}
It is an easy exercise to verify that both sides have the same divisor
and the same leading term at~$\Ocal$. Or see~\cite[Appendix~I,
Proposition~4]{MT} for a more general version.
\end{proof}

\section{Divison polynomials over $\CC$ and the Weierstrass $\s$-function}
Let~$E/\CC$ be an elliptic curve and fix an isomorphism
$\F:\CC/L\to E(\CC)$ with a lattice $L\subset\CC$. The classical
$n$-division function on~$\CC/L$ is 
\[
  \psi_n(\z,L) = \frac{\s(n\z,L)}{\s(\z,L)^{n^2}},
\]
where~$\s(\z,L)$ is the Weierstrass $\s$-function. We check that with
the given normalization of~$F_n$, the relationship between~$F_n$
and~$\psi_n$ behaves consistantly with respect to~$n$.

\begin{lemma}
\label{lemma:Fvspsi}
Let~$E/\CC$ be an elliptic curve and $\F:\CC/L\to E(\CC)$ an
isomorphism as above. Then there is a constant $\g\in\CC^*$ so that
\[
  F_n(\F(\z)) = \g^{1-n^2}\psi_n(\z,L)
  \qquad\text{for all $\z\in\CC$ and all $n\ge1$.}
\]
\end{lemma}
\begin{proof}
The analytic $n$-division function
$\psi_n(\z,L)$ on $\CC/L$ has the same divisor
as $F_n\circ\F$, so they are constant multiples of one another,
\begin{equation}
  \label{equation:Fvspsi1}
  F_n(\F(\z)) = c_n \psi_n(\z,L)\qquad\text{for all $\z\in\CC$.}
\end{equation}
The Weierstrass $\s$-function satisfies $\s(\z,L)=\z+O(\z^2)$ as
$\z\to0$, so $\psi_n(\z,L) = (n/\z^{n^2-1})\bigl(1+O(\z)\bigr)$. The
map~$\F$ has the form $(z\circ\F)(\z) = \g \z + O(\z^2)$ in a
neighborhood of~$0$, since~$\F$ is an isomorphism and~$\z$ and~$z$
are, respectively, uniformizers in neighborhoods of~$0$
and~$\Ocal=\F(0)$. Hence
\[
  F_n(\F(\z)) = F_n\bigl(\g \z+O(\z^2)\bigr) 
  = \frac{n}{(\g \z)^{n^2-1}}\bigl(1+O(\z)\bigr) 
  = \g^{1-n^2}\psi_n(\z,L).
\]
Comparing this with~\eqref{equation:Fvspsi1}, we see that 
with our chosen normalization of~$F_n$, there is
a single constant~$\g\in\CC^*$ such that
\[
  F_n\circ\F = \g^{1-n^2} \psi_n\qquad\text{for all $n\ge1$.}
  \qedtag
\]
\EndProofAtDisplay
\end{proof}

\section{Periodicity of division polynomials over finite fields}
\label{section:periodictydivpoly}

In this section we prove that the values of division polynomials over
finite fields form a purely periodic sequence. Our proof is modeled
after an analogous result by Ward (especially~\cite[Theorems~8.1
and~9.2]{Ward1}) for elliptic divisibility sequences. The proof uses a
lift to characteristic zero and the Lefschetz principle. It would
be interesting to find a purely finite field proof.

\begin{theorem}
\label{theorem:divpolyperiodic}
Let~$\FF$ be a finite field, let~$E/\FF$ be an elliptic curve, and
let~$P\in E(\FF)$ be a point of exact order~$r\ge2$. 
Then there are units~$a,b\in \FF^*$, depending on~$P$,
such that:
\begin{parts}
\Part{(a)}
If $r\ge3$, then
\begin{equation}
  \label{equation:F=powers}
  F_{kr+n}(P) =  a^{kn} b^{k^2} F_n(P) 
  \qquad\text{for all $k,n\ge0$.}
\end{equation}
\Part{(b)}
If $r=2$, then 
\begin{equation}
  \label{equation:F2kr=2}
  \left.
  \begin{aligned}
    F_{2k}(P) &= 0 \\
    F_{2k+1}(P) &= a^k b^{(k^2-k)/2} \\
  \end{aligned}
  \right\}
  \quad\text{for all $k\ge0$}
\end{equation}
\end{parts}
\end{theorem}

As an immediate corollary, we deduce the periodicity of the values of
the division polynomials.

\begin{corollary}
\label{corollary:divpolyperiodic}
Let $P\in E(\FF)$ be as in the statement of
Theorem~\ref{theorem:divpolyperiodic}. Then the sequence
\begin{equation}
  \label{equation:Fnmodp}
   \bigl(F_n(P)\bigr)_{n\ge0}
\end{equation}
is purely periodic. More precisely, if~$P$ 
has order~$r$ in~$E(\FF)$ and if we let
$q=\#\FF$, then the sequence~\eqref{equation:Fnmodp} has period~$rt$,
where~\text{$t|q-1$} if $r\ge3$ and~\text{$t|2q-2$} if $r=2$.
\end{corollary}

\begin{proof}[Proof of Corollary \ref{corollary:divpolyperiodic}]
We begin with the case that $r\ge3$.
Let~$a$ and~$b$ be as in Theorem~\ref{theorem:divpolyperiodic}, and
let~$t\ge1$ be the smallest integer such that
\[
  a^t = b^{t^2} = 1.
\]
In particular, $t$~divides the least common multiple of the orders
of~$a$ and~$b$ in~$\FF^*$, so~$t$ divides~\text{$q-1$}.
Theorem~\ref{theorem:divpolyperiodic} tells us that
\[
  F_{rt+n}(P) = a^{tn} b^{t^2} F_n(P) = F_n(P)
  \qquad\text{for all $n\ge0$,}
\]
which shows that the sequence~\eqref{equation:Fnmodp} is periodic
and that~$rt$ is a period. 
\par
Let~$\ell\ge1$ be the smallest period, i.e., the smallest integer such
that $F_{\ell+n}(P)=F_n(P)$ for all $n\ge0$. We note that $F_n(P)=0$
if and only if~$r|n$, since~$P$ has exact order~$r$ in~$E(\FF)$. From
$F_{\ell+r}(P)=F_r(P)=0$, we deduce that~$r|\ell$, say $\ell=rs$.
Since~$rt$ is a period and~$\ell=rs$ is the smallest period, we have
$s|t$, which completes the proof if $r\ge3$. (With a bit more work,
one can show that~$s=t$.)
\par
If~$r=2$, then it is easy to see from~\eqref{equation:F2kr=2}
that~$F_n(P)$ is periodic and that the period must be even.
Further we compute
\begin{align*}
  F_{2k+1+4(q-1)}(P)
  &= F_{2(k+2q-2)+1}(P) \\
  &= a^{k+2(q-1)} b^{(k^2-k)/2 + (q-1)(2k + 2(q-1) - 1)} \\
  &= a^k b^{(k^2-k)/2}\qquad\text{since $a^{q-1}=b^{q-1}=1$,} \\
  &= F_{2k+1}(P)
\end{align*}
Thus the period divides \text{$4(q-1)$} and is even, so it has the
form~$rt$ with \text{$t|2q-2$}.
\end{proof}

\begin{proof}[Proof of Theorem \ref{theorem:divpolyperiodic}]
Before starting the proof, we note that if~$r|n$, then
both~$F_{kr+n}(P)$ and~$F_n(P)$ vanish, so the desired
formula~\eqref{equation:F=powers} is vacuously true for any choice
of~$\a$ and~$\b$. We assume henceforth that~$r\notdivide n$.
\par
Let~$R$ be a complete local ring of characteristic zero with residue
field~$\FF$ (e.g., the Witt ring over~$\FF$), let~$K$ be the fraction
field of~$R$, let~$\gp$ be the maximal ideal of~$R$, and let~$\Ecal/R$
be a lift of~$E/\FF$ given by a Weierstrass equation whose reduction
modulo~$\gp$ is the Weierstrass equation of~$E/\FF$ used to normalize
the division polynomials on~$E/\FF$.
\par
The reduction map~$\Ecal(R)\to E(\FF)$ is surjective, so we can
lift~$P\in E(\FF)$ to a point in~$\Ecal(R)$. We would like to lift~$P$
to a torsion point of order~$r$. If~$r$ is not divisible by the
characteristic~$p$ of~$\FF$, then there is a unique such lift, which can
be computed by taking any lift~$Q\in\Ecal(R)$ and
computing the limit (see Proposition~\ref{proposition:teichmuller})
\[
  P' = \lim_{\substack{k\to\infty\\p^k\equiv1\MOD{r}\\}} [p^k](Q).
\]
\par
In general, if $r=p^er'$ with $p\notdivide r'$, it suffices by the
Chinese remainder theorem to lift $[p^e](P)$ and $[r'](P)$, so we are
reduced to the case that~$r$ is a power of~$p$, say $r=p^e$ with
$e\ge1$. Then it may not be possible to lift~$P$ to a torsion point
in~$\Ecal(R)$, but it is possible to do so in a finite
(ramified) extension, since we always have an exact sequence
\begin{equation}
  \label{equation:pinftorsion}
  \begin{CD}
 0 @>>> \Ecal^f(\Rbar)[p^\infty] @>>> \Ecal(\Rbar)[p^\infty]
   @>>> \Ecal(\bar\FF)[p^\infty]   @>>> 0,
  \end{CD}
\end{equation}
where~$\Ecal^f$ is the formal group of~$\Ecal$. (We also note
that~$E/\FF$ is necessarily ordinary, since~$E(\FF)$ contains the
$p^e$-torsion point~$P$.)
\par
We may thus choose a finite extension~$K'/K$ with ring of
integers~$R'/R$, residue field~$\FF'/\FF$, and maximal ideal~$\gp'|\gp$ so
that there is a torsion point $P'\in\Ecal(R')_\tors$ satisfying
$P'\cong P\MOD{\gp'}$. The point~$P'$ may not be uniquely determined
by~$P$, but this will not affect our argument.
\par
We next choose a subfield~$K'_0$ of~$K'$ that is small enough so
that we can embed~$K'_0$ into~$\CC$, but large enough so that the
given Weierstrass equation for~$\Ecal$ has coordinates in~$K'_0$ and
so that~$P'\in \Ecal(K'_0)$. Having done this, we obtain an embedding
\[
  \Ecal(K_0') \subset \Ecal(\CC) \xleftarrow[\quad\cong\quad]{\F} \CC/L
\]
for some lattice $L\subset\CC$. We let $P'=\F(\xi)$ under this
identification.
\par
Lemma~\ref{lemma:Fvspsi} tells us that there is a constant~$\g\in\CC^*$
so that
\[
  F_n(\F(\z)) = \g^{1-n^2}\frac{\s(n\z)}{\s(\z)^{n^2}}
  \qquad\text{for all $\z\in\CC$ and all $n\ge1$,}
\]
where to ease notation, we will omit reference to the lattice~$L$.
This allows us to compute the ratio of division functions as
\begin{equation}
  \label{equation:FkrnoverFn}
  F_{kr+n}(\F(\z))
  = \frac{\s(kr\z+n\z)}{\s(n\z)} \cdot \bigl(\g\s(\z)\bigr)^{n^2-(kr+n)^2}
  \cdot F_n(\F(\z)),
\end{equation}
valid for all $\z\in\CC$ with $n\z\notin L$. 
\par
By assumption, the point $\xi\in\CC/L$ has order~$r$. If we
identify~$\xi$ with a particular element of~$\CC$, then $r\xi\in
L$.  This allows us to apply the transformation formula for the
$\s$-function (see~\cite[Theorem~I.5.4]{ATAEC})
\begin{equation}
  \label{equation:sigmafncleqn}
  \s(\z+\l) = \Psi(\l)e^{\eta(\l)(\z+\l/2)}\s(\z)
  \qquad\text{for all $\z\in\CC$ and $\l\in L$.}
\end{equation}
Here~$\Psi(\l)\in\{\pm1\}$ and~$\eta(\l)$ is the quasiperiod
associated to~$\l$. More precisely,~$\Psi$ is a homomorphism
$\Psi:L/2L\to\{\pm1\}$ and~$\eta$ is a homomorphism $\eta:L\to\CC$.
Applying~\eqref{equation:sigmafncleqn} with $\z=n\xi$ and $\l=kr\xi$
and using the fact that~$\Psi$ and~$\eta$ are homomorphisms yields
\begin{align}
  \frac{\s(kr\xi+n\xi)}{\s(n\xi)}
  &= \Psi(kr\xi) e^{\eta(kr\xi)(n\xi+kr\xi/2)}  \notag\\
  &= \Psi(r\xi)^k \left(e^{\eta(r\xi)\xi}\right)^{kn}
    \left(e^{\eta(r\xi)r\xi/2}\right)^{k^2}.
  \label{equation:sigmaratio}
\end{align}
This is valid if $n\xi\notin L$, or equivalently, if~$r\notdivide n$,
since~$\xi$ has exact order~$r$ in~$\CC/L$.
\par
Now we substitute~\eqref{equation:sigmaratio}
into~\eqref{equation:FkrnoverFn} with~$\z=\xi$
to obtain
\begin{multline*}
  F_{kr+n}(\F(\xi)) \\
  = \Psi(r\xi)^k \left(e^{\eta(r\xi)\xi}(\g\s(\xi))^{-2r}\right)^{kn}
    \left(e^{\eta(r\xi)r\xi/2}(\g\s(\xi))^{-r^2}\right)^{k^2} F_n(\F(\xi)).
\end{multline*}
In other words, we have proven that there exist numbers~$\a,\b\in\CC$,
depending only on~$\xi$ and independent of~$k$ and~$n$, so that
\[
  F_{kr+n}(\F(\xi)) = \a^{kn} \b^{k^2}F_n(\F(\xi))
  \qquad\text{for all $k,n\ge0$ with $r\notdivide n$.}
\]
(We have absorbed $\Psi(r\xi)^k=(\pm1)^k$ into the $\b^{k^2}$ terms.)
\par
Recall that~$\xi\in\CC/L$ corresponds to the point
$P'=\F(\xi)\in\Ecal(K_0')$, so we may equally well write this as
\begin{equation}
  \label{equation:FTformula}
  F_{kr+n}(P') = \a^{kn} \b^{k^2}F_n(P')
  \qquad\text{for all $k,n\ge0$.}
\end{equation}
(We drop the restriction that $r\notdivide n$, since as noted earlier,
the formula~\eqref{equation:FTformula} is trivially true in this
case.)  
\par
We now make the assumption that $r\ge3$, and at the end we will
briefly indicate the changes needed to deal with the case $r=2$.
We substitute $(k,n)=(1,1)$ and $(k,n)=(1,2)$ into~\eqref{equation:FTformula}
to obtain
\[
  F_{r+1}(P') = \a\b F_1(P')=\a\b
  \qquad\text{and}\qquad
  F_{r+2}(P') = \a^2\b F_2(P').
\]
(Note that~$F_1=1$.) Our assumption that $r\ge3$ implies
that~$F_2(P')\ne0$, so we can solve for~$\a$ and~$\b$,
\begin{equation}
  \label{equation:alphabeta}
  \begin{aligned}
  \a&=\frac{F_{r+2}(P')}{F_2(P')F_{r+1}(P')}\in K(P')\subset K',\\
  \b&=\frac{F_2(P')F_{r+1}(P')^2}{F_{r+2}(P')}\in K(P')\subset K'.\\
  \end{aligned}
\end{equation}
Thus we may view~\eqref{equation:FTformula} as a formula in the
complete local field~$K'$, since all of the quanitities appearing in
it are in~$K'$.
\par
We claim that~$\a$ and~$\b$ are actually $\gp'$-units in~$K'$. This
follows from the fact that for points
$Q\in\Ecal(R')\setminus\Ecal^f(R')$ and integers $n\ge1$, we have
\[
  F_n(Q)\equiv 0\pmod{\gp'}
  \qquad\text{if and only if}\qquad
  nQ\equiv \Ocal\pmod{\gp'}.
\]
Thus~\eqref{equation:alphabeta} shows that~$\a$ and~$\b$ are
$\gp'$-units provided that
\begin{equation}
  \label{equation:abarepunits}
  \begin{aligned}
  (r+2)P'&\not\equiv\Ocal\pmod{\gp'},\\
  (r+1)P'&\not\equiv\Ocal\pmod{\gp'},\\ 
  2P'&\not\equiv\Ocal\pmod{\gp'}.\\
  \end{aligned}
\end{equation}
But~$P'$ modulo~$\gp'$ has  exact period $r\ge3$,
so the three conditions~\eqref{equation:abarepunits} are true.
\par
To recapitulate, we have shown that there are $\gp'$-units $\a,\b\in K'$
such that
\begin{equation}
  \label{equation:FTformulainL}
  F_{kr+n}(P')
  =\a^{kn} \b^{k^2}F_n(P')
  \qquad\text{for all $k,n\ge0$.}
\end{equation}
We reduce this formula modulo~$\gp'$ and use the fact that $P'\equiv
P\MOD{\gp}$ (remember that we chose~$P'$ to be a lift of the
point~$P\in E(\FF)$) to obtain
\[
  F_{kr+n}(P)
  = a^{kn} b^{k^2} F_n(P)
  \qquad\text{for all $k,n\ge0$,}
\]
where~$a$ and~$b$ are elements of the residue field~$\FF'$ of~$K'$. To
see that~$a,b\in\FF$, we substitute $(k,n)=(1,1)$ and $(k,n)=(1,2)$
and solve for~$a,b$ (cf.~\eqref{equation:alphabeta}) to obtain
\begin{equation}
  \label{equation:ab}
  a=\frac{F_{r+2}(P)}{F_2(P)F_{r+1}(P)}\in\FF,\quad
  b=\frac{F_2(P)F_{r+1}(P)^2}{F_{r+2}(P)}\in\FF.
\end{equation}
(The~$F_i(P)$ values are nonzero, since $F_i(P)\equiv
F_i(P')\MOD{\gp'}$.)  This completes the proof of
Theorem~\ref{theorem:divpolyperiodic} for $r\ge3$.
\par
Suppose now that $r=2$. Then it is not helpful to substitute $(k,n)=(1,2)$
into~\eqref{equation:FTformula}, since both sides are zero. So instead we 
substitute~$(k,n)=(1,1)$ and~$(k,n)=(1,3)$ to obatin
\[
  F_{r+1}(P') = \a\b F_1(P')=\a\b
  \qquad\text{and}\qquad
  F_{r+3}(P') = \a^3\b F_3(P'),
\]
where we know that~$F_3(P')\ne0$ since~$P'$ has order~$r=2$.
We can no longer solve for~$\a$ and~$\b$, but we can solve for
for the quantities~(cf.~\eqref{equation:alphabeta})
\begin{equation}
  \b^2=\frac{F_3(P')F_{r+1}(P')^3}{F_{r+3}(P')}
  \quad\text{and}\quad
  \a\b=F_{r+1}(P').
\end{equation}
Thus~$\b^2$ and~$\a\b$ are in~$K'$, and by the same argument given
earlier, they are actually~$\gp'$ units in~$K'$.  We set~$n=1$
and~$r=2$ in~\eqref{equation:FTformulainL} to obtain
\[
  F_{2k+1}(P') =\a^k \b^{k^2} = (\a\b)^k (\b^2)^{(k^2-k)/2}.
\]
Reducing this formula modulo~$\gp'$ and using the fact that~$P'\bmod\gp'$
is equal to~$P$, we see that there are 
units $a,b\in\FF'$ (i.e., 
\text{$a\equiv\a\b\MOD{\gp'}$}
and 
\text{$b\equiv\b^2\MOD{\gp'}$}) so that
\[
  F_{2k+1}(P) =  a^k b^{(k^2-k)/2}
  \qquad\text{for all $k\ge0$.}
\]
Finally, putting $k=1$ shows that $a=F_3(P)\in\FF$ and then putting
$k=2$ shows that $b=a^{-1}F_5(P)\in\FF$, which completes the
proof of Theorem~\ref{theorem:divpolyperiodic} for $r=3$.
\end{proof}

\section{The Teichm\"uller character}
\label{section:teichmuller}

The classical \emph{Teichim\"uller character} is the unique homomorphism
\[
  \chi : \ZZ_p^*\longrightarrow\bfmu_{p-1}
  \quad\text{satsifying}\quad
  \chi(a)\cong a \pmod{p}.
\]
The Teichm\"uller character may be constructed as $\chi(a)=\lim
a^{p^k}$. It is well known how to generalize this construction to
group schemes~$G$ over~$\ZZ_p$ or other complete local rings.

\begin{proposition}
\label{proposition:teichmuller}
Let $K/\QQ_p$ be a finite extension, let~$R$ be the ring of integers
of~$K$, let~$\gp$ be the maximal ideal of~$R$, and let~$\FF$ be the
residue field of~$R$.  Let~$G/R$ be a group scheme, and for
any point~$a\in G(R)$, let~$\t(a)$ denote
the order of~$a\bmod p$ in the special fiber~$G(\FF)$.
We denote by
\[
  G'(R) = \bigl\{ a\in G(R) : p\notdivide \t(a) \bigr\}
\]
the pullback to~$G(R)$ of the prime-to-$p$ part of~$G(\FF)$.
\begin{parts}
\Part{(a)}
There is a unique homomorphism
\begin{equation}
  \label{equation:teichprops}
  \chi : G'(R)\longrightarrow G(R)_\tors
  \quad\text{satsifying}\quad
  \chi(a)\equiv a \pmod{\gp}.
\end{equation}
We call~$\chi$ the \emph{Teichm\"uller ``character''} for~$G/R$.
\Part{(b)}
Writing the gruop law in~$G(R)$ multiplicatively, 
the \emph{Teichm\"uller character} can be computed as the limit
\begin{equation}
  \label{equation:teichlimit}
  \chi(a) = \lim_{\substack{k\to\infty\\p^k\equiv1\MOD{\t(a)}\\}} a^{p^k}.
\end{equation}
\Part{(c)}
The order of~$\chi(a)$ in $G(R)_\tors$ is exactly~$\t(a)$.
\Part{(d)}
The reduction map $G'(R)_\tors\to G'(\FF)$ is an isomorphism.
\end{parts}  
\end{proposition}
\begin{proof}
For lack of a suitable reference, we sketch the short proof of this
well-known result. We begin by proving that the
limit~\eqref{equation:teichlimit} in~(b) exists.  
\par
Let $G_1(R)$ be the kernel of the reduction map $G(R)\to
G(\FF)$. Then $G_1(R)$ is a pro-$p$ group.  Let~$q=p^e$ be the
smallest power of~$p$ satisfying $q\equiv1\MOD{\t(a)}$. Then
$a^q\equiv a\MOD{\gp}$, so for $i>j$ we have
\[
  a^{q^i}\cdot a^{-q^j} = \bigl(a^{q^{i-j}-1}\bigr)^{q^j} \longrightarrow 1_G
  \qquad\text{as $i>j\to\infty$.}
\]
This is true since $a^{q^{i-j}}\equiv a\MOD{\gp}$, so
$a^{q^{i-j}-1}\in G_1(R)$. Thus the
sequence~\eqref{equation:teichlimit} is Cauchy, so it converges.
\par
For~$a\in G'(R)$, we now define~$\chi(a)$ to be the
limit~\eqref{equation:teichlimit}.  Then
\[
  \chi(a)=\lim_{i\to\infty} a^{q^i}\equiv a\MOD{p}
  \qquad\text{and}\qquad
  \chi(a)^q=\chi(a).
\]
In particular, $\chi(a)^{q-1}=1_G$, so~$\chi(a)\in G(R)_\tors$.
This shows that the function~$\chi$ defined by the
limit~\eqref{equation:teichlimit} has the desired
properties~\eqref{equation:teichprops}, and it is obvious that it is a
homomorphism, which proves the existence part of~(a). The uniqueness
is immediate from the fact that~$G_1(R)$ has no prime-to-$p$
torsion and that~$\chi(a)$ has order dividing~\text{$q-1$}. 
\par
Let~$T$ be the order of~$\chi(a)$, so from above, 
\text{$T|q-1$}. In particular,~$T$ is prime to~$p$. Further, we have
\[
  \chi(a)^{\t(a)} \equiv a^{\t(a)}  \equiv 1 \pmod{\gp},
\]
so $\chi(a)^{\t(a)}$ is a $T^{\text{th}}$~root of unity whose reduction
modulo~$\gp$ is~$1$. The formal multiplicative group~$\GG_{m,1}(R)$
has no prime-to-$p$ torsion, so $\chi(a)^{\t(a)}=1$. Thus $T|\t(a)$.
Conversely, $1=\chi(a)^T\equiv a^T\equiv1\MOD{\gp}$, so~$\t(a)|T$. 
Hence $T=\t(a)$, which completes the proof of~(c).
\par
The proof of~(d) is also immediate. The map~$G'(R)_\tors\to
G(\FF)$ is injective, since the formal group~$G_1(R)$ has no
prime-to-$p$ torsion. On the other hand, let~$\a\in G(\FF)$ have
order~$\t$ with $p\notdivide\t$. The group~$G$ is smooth over~$R$,
so we can choose an $a\in G(R)$ with $a\equiv\a\MOD{p}$. Then
$\t(a)=\t$ by definition, and $\chi(a)\in G(R)_\tors$ satisfies
\[
  \chi(a)\equiv a\equiv\a\MOD{\gp}.
\]
This shows that the map $G'(R)_\tors\to G'(\FF)$ is surjective, which
completes the proof of~(d).
\end{proof}

\section{The Mazur-Tate $p$-adic sigma function}
\label{section:MTsigmafunction}

In this section we recall a construction of Mazur and Tate, and in the
next section we apply their construction to prove a $p$-adic limit for division
polynomials on elliptic curves with ordinary reduction. We set the
following notation (following~\cite{MT}):

\begin{notation}
\item[$K$]
a finite extension of~$\QQ_p$.
\item[$R$] 
the ring of integers of~$K$.
\item[$\gp$]
the maximal ideal of~$R$.
\item[$\FF$]
the residue field $R/\gp$ of $K$.
\item[$\Rbar$]
the integral closure of~$R$ in an algebraic closure~$\Kbar$ of~$K$.
\item[$E/K$]
an elliptic curve over $K$. We also fix a minimal Weierstrass equation
for~$E/K$, from which we obtain an invariant differential
$\o=dx/(2y+a_1x+a_3)$ and a uniformizer $z=-x/y$ at~$\Ocal$ satisfying
$(\o/dz)(\Ocal)=1$.
\item[$\Ecal/R$]
the N\'eron model of~$E$.
\item[$\Ecal^f$]
the formal group of~$\Ecal$.
\item[$F_n$]
the $n$-division polynomial on~$E$, that is, the rational
function~$F_n\in K(E)$ satisfying Definition~\ref{definition:ndivpoly}
for the given Weierstrass equation.
\end{notation}

\begin{remark}
\label{remark:MTdivpolydef}
Mazur and Tate~\cite{MT} define division polynomials~$F_n\in K(E)$
to be the unique functions satisfying
\[
  (F_n) = [n]^*(\Ocal) - n^2(\Ocal)
  \quad\text{and}\quad
  \left(\frac{z^{n^2}F_n(z)}{[n](z)}\right)(\Ocal)=
  \left(\frac{\o}{dz}(\Ocal)\right)^{1-n^2}.
\]
This agrees with our definition~\eqref{definition:ndivpoly}, since we
have chosen~$z$ and~$\o$ compatibly to satisfy $(\o/dz)(\Ocal)=1$.
\end{remark}

\begin{theorem}[Mazur-Tate]
\label{theorem:MTsigma}
With the above notation and normalizations, assume that $p\ge3$ and
that~$\Ecal$ has good ordinary reduction.  Then there is a unique
power series $\s\in z+z^2R[[z]]$ satisfying
\begin{equation}
  \label{equation:MTsigma}
  \s(nQ) = \s(Q)^{n^2} F_n(Q)
  \qquad\text{for all $Q\in \Ecal^f(\Rbar)$.}
\end{equation}
\end{theorem}
\begin{proof}
See \cite[Section~2]{MT} for the construction of~$\s$ and
\cite[Theorem~3.1]{MT} for a description of its properties.  The
construction of~$\s$ actually works as long as~$\Ecal$ has ordinary
reduction, i.e., if $\Ecal^f$ is isomorphic over~$\bar\FF$ to the formal
multiplicative group~$\GG_m^f$. If one chooses a different Weierstrass
equation for~$E$, then~$F_n$ changes by a constant factor of the
form~$c^{n^2-1}$ with $c\in R^*$, and hence~$\s$ changes by a
factor~$c$.
\end{proof}

\begin{remark}
Theorem~\ref{theorem:MTsigma} remains true for~$p=2$ provided that
everything is squared. That is, there is a unique power series
$\s^2\in z^2+z^3R[[z]]$ satisfying $\s^2(nQ) = \s(Q)^{2n^2}
F_n^2(Q)$. But it is not possible to unambiguously take a square root
and have~\eqref{equation:MTsigma} hold for all~$n\ge1$ and
all~$Q\in\Ecal^f(\Rbar)$.
\end{remark}

\section{A $p$-adic limit of division polynomials}
\label{section:padiclimitofdivpoly}

In this section we compute the $p$-adic limit of the values of certain
subseqeunces of the division polynomials evaluated at a point. We
continue with the notation from Section~\ref{section:MTsigmafunction}.

\begin{theorem}
\label{theorem:divpolyconverge}
Assume that~$p\ge3$ and that~$\Ecal$ has good ordinary reduction, and
let $P\in\Ecal(R)\setminus\Ecal^f(R)$.  
\begin{parts}
\Part{(a)}
There exists a power~$q=p^N$ so that for every $m\ge1$, the limit
\begin{equation}
  \label{equation:GmqP}
  G_{m,q}(P) := \lim_{k\to\infty} F_{mq^k}(P)\quad\text{converges in $R$.}
\end{equation}
\Part{(b)}
$G_{m,q}(P)$ is algebraic over~$\QQ(E)$. 
\Part{(c)}
In order to specify an allowable value for~$q$, let \text{$r\ge2$} be
the order of~$P\bmod\gp$. Then
Corollary~\ref{corollary:divpolyperiodic} tells us that the sequence
\[
  (F_n(P)\bmod\gp)_{n\ge0}
\]
is periodic with period~$rt$ for some integer~\text{$t\ge1$}
with~$p\notdivide t$.  Let~$r'$ be the $p$-free part of~$r$, that is,
$r'=rp^{-\ord_p(r)}$, and let~$e\ge1$ be an exponent so that
\[
  q=(\Norm\gp)^e\quad\text{satisfies}\quad
  q\equiv1\MOD{r't}.
\]
Then the limit~\eqref{equation:GmqP} in~{\upshape(}a{\upshape)}
exists for this value of~$q$.
\Part{(d)}
Continuing with the notation from~{\upshape(}c{\upshape)}, we have
\[
  G_{m,q}(P)^{{r'}^2} \in \QQ(\bfmu(K),E[r']),
\]
where~$\bfmu(K)$ denotes the roots of unity in~$K$.  Further,
\[
  G_{m,q}(P)=0\quad\text{if and only if}\quad m\equiv0\MOD{r'}.
\]
\end{parts}
\end{theorem}

\begin{remark}
We note that it is quite easy to estimate the valuation of~$F_n(P)$,
either directly as in~\cite{CH} or using the transformation formula
for local height functions. In particular, with notation as in
Theorem~\ref{theorem:divpolyconverge}, it is an elementary
exercise to prove that
\[
  v\bigl(F_{rn}(P)\bigr) = v(n) + O(1)\quad\text{for all $n\ge1$,}
\]
and hence~$F_{rp^k}(P)\to 0$ as $k\to\infty$.  Thus the interest and
the depth of Theorem~\ref{theorem:divpolyconverge} lies in the
convergence of~$F_{mq^k}(P)$ in those cases that the limit is not
zero.
\end{remark}

Before starting the proof of Theorem~\ref{theorem:divpolyconverge}, we
give an elementary result that will allow us to take roots of
convergent sequences.

\begin{lemma}
\label{lemma:seqrootconv}
Let $(A_k)_{k\ge0}$ be a sequence in~$R^*$ with the property
that $(A_k\bmod\gp)_{k\ge0}$ is constant. Let~$n\ge1$ be an integer
with~$p\notdivide n$. Then
\[
  \lim_{k\to\infty} A_k^n~\text{exists in $R$}
  \Longleftrightarrow
  \lim_{k\to\infty} A_k~\text{exists in $R$.}
\]
\end{lemma}
\begin{proof}
One direction is trivial. So we assume that $\lim_{k\to\infty}A_k^n$
exists and we must prove that we may take the $n^{\text{th}}$~root.
Fix~$\a\in R^*$ with \text{$\a\equiv A_k\MOD{\gp}$} for all $k\ge0$. Then
for any $j,k\ge0$ we have
\[
  \sum_{\ell=0}^{n-1} A_{j}^\ell\cdot A_{k}^{n-1-\ell}
  \equiv n\a^{n-1} \not\equiv 0 \pmod{\gp},
\]
from which we deduce that
\begin{align*}
  \lim_{j,k\to\infty} \left|A_{j}-A_{k}\right|
  &= \lim_{j,k\to\infty}\frac{\left|A_{j}^n-A_{k}^n\right|}
   {\displaystyle\biggl|\sum_{\ell=0}^{n-1} 
     A_{j}^\ell\cdot A_{k}^{n-1-\ell}\biggr|}\\
  &= \lim_{j,k\to\infty}\left|A_{j}^n-A_{k}^n\right|\\
  &=0\qquad\text{since $A_k^n$ converges as $k\to\infty$.}
\end{align*}
This shows that the sequence~$(A_{k})_{k\ge0}$ is Cauchy, hence converges,
in~$R$, which completes the proof of the lemma.
\end{proof}

We are now ready to prove our main result.

\begin{proof}[Proof of Theorem \ref{theorem:divpolyconverge}]
We use Lemma~\ref{lemma:divpolychainrule} twice to obtain
\[
 \bigl(F_r\circ[n]\bigr)\cdot F_n^{r^2} =
  F_{nr} = F_{rn} = \bigl(F_n\circ[r]\bigr)\cdot F_r^{n^2},
\]
and hence
\[
  F_n^{r^2} 
  = \frac{\bigl(F_n\circ[r]\bigr)\cdot F_r^{n^2}}{\bigl(F_r\circ[n]\bigr)}.
\]
We evaluate this identity at the point~$P$ and use the fact
that $rP\in\Ecal^f(R)$ to rewrite~$F_n(rP)$
using the Mazur-Tate sigma function. Thus
\begin{align}
  \label{equation:fnPr2}
  F_n(P)^{r^2}
  &= \frac{F_n(rP)F_r(P)^{n^2}}{F_r(nP)} \notag\\
  &= \frac{\left(\s(nrP)/\s(rP)^{n^2}\right)F_r(P)^{n^2}}{F_r(nP)} 
         \qquad\text{from Theorem~\ref{theorem:MTsigma},} \notag\\
  &= \frac{\s(rnP)}{F_r(nP)}\cdot
    \left(\frac{F_r(P)}{\s(rP)}\right)^{n^2}.
\end{align}
\par
We consider first the case that~$p\notdivide r$, so we may let
$T=\chi(P)\in E[r]$ be the Teichm\"uller image of~$P$
(Proposition~\ref{proposition:teichmuller}), and then $Q=P-T$
satisfies $Q\in\Ecal^f(R)$.  (See Section~\ref{section:teichmuller}.)
Note that $T\ne\Ocal$, since we have assumed that~$P\notin\Ecal^f$.
Using $rP=r(T+Q)=rQ$, we can rewrite~\eqref{equation:fnPr2} as
\begin{equation}
  \label{equation:fnPr2v2}
  F_n(P)^{r^2}
   = \frac{\s(rnQ)}{F_r(nT+nQ)}\cdot
    \left(\frac{F_r(T+Q)}{\s(rQ)}\right)^{n^2}.
\end{equation}
\par
Let
\[
  \t_T:\Ecal\longrightarrow\Ecal
\]
be the translation-by-$T$ map. The division function~$F_r$ has simple
zeros at all nonzero $r$-torsion points, and our assumption
that~$p\notdivide r$ implies that the same is true of the restriction
of~$F_r$ to the special fiber of~$\Ecal$. Hence
\begin{equation}
  \label{equation:Ftau=gzT}
  F_r\circ\t_T = z \cdot g_T
\end{equation}
for a rational function~$g_T\in K(E)$ whose restriction to the special
fiber of~$\Ecal$ is regular and nonvanishing at~$\Ocal$, i.e.,
$g_T(\Ocal)\in R^*$.
\par
Then using the fact that $Q\in\Ecal^f(R)$ and $p\notdivide r$, we see
that
\begin{align*}
  \frac{F_r(T+Q)}{\s(rQ)} &= \frac{z(Q)\cdot g_T(Q)}{\s(rQ)} \\
  &= \frac{z(Q)}{z(rQ)} \cdot \frac{z(rQ)}{\s(rQ)} \cdot g_T(Q) \\
  &\equiv \frac{1}{r} g_T(\Ocal) \pmod{\gp}.
\end{align*}
Hence if we let~$n=mq^k$ with~$q$ a power of~$\Norm\gp$
and with some fixed~$m$ and if we let~$k\to\infty$, then we
can evaluate the limit of the second factor
in~\eqref{equation:fnPr2v2} as
\[
  \lim_{k\to\infty} \left(\frac{F_r(T+Q)}{\s(rQ)}\right)^{(mq^{k})^2}
  = \chi\left(\frac{g_T(\Ocal)}{r}\right)^{m^2},
\]
where $\chi:R^*\to\bfmu(K)$ is the Teichm\"uller character on~$K$
(cf{.} Proposition~\ref{proposition:teichmuller}). In particular,
the value is a root of unity in~$K$.
\par
In order to evaluate the limit of the first factor
in~\eqref{equation:fnPr2v2}, we take a sequence of~$n$'s
of the form $n=mq^k$ with $k=1,2,3,\dots$, 
where~$q$ is a certain fixed power of~$p$. More precisely,
we already noted that we want~$q$ to be a power of~$\Norm\gp$,
and we now further specify that
\begin{equation}
  \label{equation:qdef}
  q=(\Norm\gp)^e\qquad\text{satisfies}\qquad
  q\equiv1\pmod{rt}.
\end{equation}
In particular, $q\equiv1\MOD{r}$, so $n\equiv m\MOD{r}$ for all~$k$,
and hence $nT=mT$ is independent of~$k$. To ease notation, we let
$Q_k=nQ=mq^kQ$, so the first factor in~\eqref{equation:fnPr2v2} is
$\s(rQ_k)/F_r(mT+Q_k)$. Our task is to evaluate the limit of this
fraction as $Q_k\to\Ocal$. There are two cases to consider.
\par  
First, if~$r|m$, then $mT=\Ocal$, so we must evaluate the limit
of $\s(rQ_k)/F_r(Q_k)$. The function~$\s\circ[r]$ has a simple zero
at~$\Ocal$, while~$F_r$ has a pole of order~$r^2-1$ at~$\Ocal$, so 
$(\s\circ[r])/F_r$ vanishes (to order~$r^2$) at~$\Ocal$. Hence in
this case the limit is~$0$.
\par
The more interesting case is when~$r\notdivide m$, so~$mT$ is a
nonzero torsion point. Then
\begin{align*}
  \frac{\s(rQ_k)}{F_r(mT+Q_k)}
  &=\frac{\s(rQ_k)}{z(Q_k)\cdot g_{mT}(Q_k)} 
     \qquad\text{from \eqref{equation:Ftau=gzT},} \\
  &=\frac{\s(rQ_k)}{z(rQ_k)} \cdot\frac{z(rQ_k)}{z(Q_k)} \cdot 
          \frac{1}{g_{mT}(Q_k)} \\
  &\xrightarrow[\;Q_k\to\Ocal\;]{} \frac{r}{g_{mT}(\Ocal)}.
\end{align*}
\par
To recapitulate, taking~$q$ as specified in~\eqref{equation:qdef},
we have proven that
\begin{equation}
  \label{equation:limitr2}
  \lim_{k\to\infty} F_{mq^{k}}(P)^{r^2} =
  \begin{cases}
    0&\text{if $r|m$,} \\
    (r/g_{mT}(\Ocal))\cdot \chi(g_T(\Ocal)/r)^{m^2}\ne 0
     &\text{if $r\notdivide m$.} \\
  \end{cases}
\end{equation}
This almost completes the proof when~$p\notdivide r$, the only
difficulty being that if~$r\notdivide m$, then we have only computed a
power of the desired limit. In order to take the $r^2$-root, we
consider the sequence
\[
  A_k = F_{mq^k}(P),\qquad k=0,1,2,\ldots.
\]
We observe first that the sequence~$A_k\bmod\gp$ is actually constant,
since the sequence~$(F_n(P)\bmod\gp)_{n\ge0}$ is periodic
(Corollary~\ref{corollary:divpolyperiodic}) with period~$rt$ and~$q$
satisfies \text{$q\equiv1\MOD{rt}$} from~\eqref{equation:qdef}.
Thus
\[
  mq^k\equiv m\MOD{rt}\qquad\text{for every~$k\ge0$,} 
\]
and hence
\[
  A_k = F_{mq^k}(P) \equiv F_m(P) = A_0 \pmod{\gp}
  \qquad\text{for all $k\ge0$.}
\] 
On the other hand, we have already proven that 
$\lim_{k\to\infty}A_k^{r^2}$ exists in~$R$.  So we can apply
Lemma~\ref{lemma:seqrootconv} (with~$n=r^2$) to the sequence~$A_k$ and
deduce that~$\lim_{k\to\infty} A_k$ exists in~$R$. This completes the
proof of the theorem in the case that~$p\notdivide r$.
\par
Next we consider the case that \text{$r=p^j$} is a power of~$p$, so in
particular~$j\ge1$.  For each integer $k\ge0$, let
$K_k=K(E[p^{k+1}])$, let~$R_k$ be the ring of integers of~$K_k$,
let~$\gp_k|\gp$ be the maximal ideal of~$R_k$, and
let~$\FF_k=R_k/\gp_k$ be the residue field of~$K_k$.
\par
We note that for any integers~$m|n$, the quotient $F_n/F_m$ is 
regular away from~$\Ocal$, since its divisor is
\[
  (F_n/F_m) = [m]^*\bigl([n/m]^*(\Ocal)-(\Ocal)) - (n^2-m^2)(\Ocal)
  \ge - (n^2-m^2)(\Ocal).
\]
In particular,
\begin{equation}
  \label{equation:Fpk+1=Fpkfpk}
  F_{p^{k+1}} = F_{p^k}\cdot f_{p^k}
\end{equation}
for a function~$f_{p^k}\in K(E)$ that is regular away from~$\Ocal$ and
that vanishes at the points in~\text{$E[p^{k+1}]\setminus E[p^k]$}.
\par
We claim that $F_{p^k}(P)\to 0$ as $k\to\infty$. Fix a point~$T\in
E[p^j]$ satisfying $T\equiv P\pmod{\gp_j}$ (cf. the exact
sequence~\eqref{equation:pinftorsion}) and let $k\ge j$.  Then all of
the points in the set
\[
  \bigl\{  T + T' : T' \in E^f[p^{k+1}] \setminus E^f[p^k]  \bigr\}
\]
have that property that
\[
  T+T' \equiv T \equiv P \pmod{\gp_k}
  \qquad\text{and}\qquad
  T+T' \in E[p^{k+1}] \setminus E[p^k].
\]
It follows that $f_k(T+T')=0$, and hence that
\[
  f_k(P) \equiv f_k(T+T') = 0 \pmod{\gp_k}.
\]
However, $f_k(P)\in K$, so $f_k(P)\equiv0\MOD{\gp}$.
Now we use~\eqref{equation:Fpk+1=Fpkfpk} repeatedly
to deduce that
\[
  F_{p^{k}}(P) = F_{p^j}(P) \prod_{i=j}^{k-1} f_{p^i}(P)
  \equiv 0 \pmod{\gp^{k-j}}
  \qquad\text{for all $k\ge j$.}
\]
Hence
\begin{equation}
  \label{equation:Fpk(P)=0}
  \lim_{k\to\infty} F_{p^k}(P) = 0.
\end{equation}
\par
More generally, still assuming that $r=p^j$, let~$m\ge1$ be a fixed
integer and write $m=m'p^\ell$ with $p\notdivide m'$. Then
\begin{alignat*}{2}
  \lim_{k\to\infty} F_{mp^k}(P)
  &= \lim_{k\to\infty} F_{m'p^k}(P) \\
  &= \lim_{k\to\infty} F_{p^k}(m'P) F_{m'}(P)^{p^k} 
    &\qquad&\text{from Lemma~\ref{lemma:divpolychainrule},} \\
  &= 0 &&\text{from \eqref{equation:Fpk(P)=0},}
\end{alignat*}
since~$m'P$ also has order~$p^j$ and $F_{m'}(P)\in R$. This completes
the proof in the case that~$r=p^j$ is a power of~$p$.
\par
Finally, we consider the case that~$p|r$, but~$r$ is not a power
of~$p$, say $r=p^jr'$ with $j\ge1$ and $r'\ge2$.  Notice that the
point $P'=p^jP$ has exact order~$r'$ modulo~$\gp$, where~$r'\ge2$
and~$p\notdivide r'$. From above, there is a power~$q$ of~$p$
so that
\[
  \lim_{k\to\infty} F_{mq^k}(P')\quad\text{exists},
\]
and further the limit is~$0$ if and only if~$r'|m$. We next use
Lemma~\ref{lemma:divpolychainrule} to write
\[
   F_{mq^k}(P)
  =  F_{mp^{-j}q^k}(p^jP) F_{p^j}(P)^{(mp^{-j}q^k)^2}
  =  F_{mp^{-j}q^k}(P') F_{p^j}(P)^{(mp^{-j}q^k)^2}.  
\]
Since we are going to let $k\to\infty$, we can pull off some powers
of~$q$ to compensate for the~$p^{-j}$. To simplify notation, fix an
exponent~$\ell$ so that~$p^j|q^\ell$ and let~$m'=mq^\ell$. Then we
find that
\begin{equation}
  \label{equation:r=pjr}
  \lim_{k\to\infty} F_{mq^k}(P)
  = \lim_{k\to\infty} F_{m'q^k}(P') F_{p^j}(P)^{(m'q^k)^2}.
\end{equation}
The point~$P'$ has exact order~$r'$ modulo~$\gp$, where~$r'\ge2$
and~$p\notdivide r'$, so from above, we know that the first
term~$F_{m'q^k}(P')$ in~\eqref{equation:r=pjr} has a limit in~$\ZZ_p$
as~$k\to\infty$, and further that the limit is~$0$ if and only
if~$r'|m$.  (Note that $r'|m$ if and only if $r'|m'$, since
$p\notdivide r'$.)  Similarly, the second term
in~\eqref{equation:r=pjr} has a limit in~$R$, since~$F_{p^j}(P)\in
R^*$ (where we again use the assumption that~$r'\ge2$). More
precisely, the limit is a root of unity, a power of the value of the
Teichm\"uller character~$\chi(F_{p^j}(P))$. Hence the limit
in~\eqref{equation:r=pjr} exists, which completes the proof of
Theorem~\ref{theorem:divpolyconverge}(a) that in all
cases,
\[
  G_{m,q}(P) = \lim_{k\to\infty}F_{mq^k}(P)
  \quad\text{exists in $R$.}
\]
\par
However, a closer examination of the proof given above shows that we
have actually completed the proof of all four parts of
Theorem~\ref{theorem:divpolyconverge}. We showed that the limit exists
for the value of~$q$ specified in~(c), and we showed
that~$G_{m,q}(P)=0$ precisely as specified in~(d). Further the value
of the limit~$G_{m,q}(P)^{{r'}^2}$ is given explicitly in terms of
certain rational functions in~$K(E)$ evaluated at certain points
in~$E[r']$, together with certain roots of unity,
so~$G_{m,q}(P)$ is algebraic over~$\QQ(E)$, and in fact satisfies
the property described in~(d).
\end{proof}

\section{Periodicity of division polynomials modulo $\gp^\mu$}

We continue with the notation used in
Sections~\ref{section:MTsigmafunction}
and~\ref{section:padiclimitofdivpoly}, so~$K/\QQ_p$ is a finite
extension and~$E/K$ an elliptic curve.  For simplicity, assume that
$p\ge3$ and that~$E$ has good reduction.  

Let~$P\in E(K)$ be a point whose reduction modulo~$\gp$ has order
$r\ge2$. We proved in Corollary~\ref{corollary:divpolyperiodic} that
the sequence $(F_n(P)\bmod\gp)$ is periodic with period~$rt$, where
$\gcd(p,t)=1$. More precisely, in
Theorem~\ref{theorem:divpolyperiodic} we gave an explicit formula
for~$F_{kr+n}(P)\bmod\gp$ as a function of~$k$ and~$n$. Of course,
when~$n=0$, then~$F_{kr}(P)\equiv0\MOD{\gp}$.

In the context of elliptic divisibility sequences, which we will study
in Section~\ref{section:EDSpadiclimit},
Shipsey~\cite[Theorem~3.5.4]{SHIP} gives a formula (when $K=\QQ_p$)
for the value of~$F_{kr}(P)$ modulo~$p^2$, and from this she
immediately deduces the periodicity of the sequence
\text{$(F_{kr}(P)\bmod p^2)_{k\ge1}$}.  We will use the Mazur-Tate
$\s$-function to prove a result that is both much stronger, and yet
not as general, as that of Shipsey. More precisely, we will prove the
periodicity of $(F_{kr}(P)\bmod\gp^\mu)_{k\ge1}$ for every fixed prime
power~$\gp^\mu$, but our proof will only be valid when~$E$ has good
ordinary reduction.

\begin{theorem}
\label{theorem:Fnperiodmodptothemu}
With notation as in Section~\ref{section:MTsigmafunction}, assume that
$p\ge3$ and that~$\Ecal$ has good ordinary reduction, and let
$P\in\Ecal(R)\setminus\Ecal^f(R)$. Further let $r\ge2$ be the order of
$P\bmod\gp$. Then for any exponent $\mu\ge1$, the sequence
\begin{equation}
  \label{equation:Fnperiodmodptothemu}
  (F_{kr}(P)\bmod\gp^\mu)_{k\ge1}\quad\text{is periodic.}
\end{equation}
More precisely, let $e=\ord_\gp(p)$ be the ramification index
of~$K_\gp/\QQ_p$ and let~$\l$ be the smallest positive integer
satisfying
\begin{equation}
  \label{equation:deflambda}
  \min_{0\le i\le \l} \bigl\{ (\l-i)e+p^i \bigr\} \ge \mu.
\end{equation}
Then the sequence~\eqref{equation:Fnperiodmodptothemu} has
period dividing $(\Norm\gp-1)p^\l$.
\end{theorem}
\begin{proof}
The point~$rP$ is in the formal group~$\Ecal^f(R)$, so the
Mazur-Tate $\s$-function can be evaluated at~$rP$. This allows us to
compute
\begin{alignat}{2}
  \label{equation:Fkr(P)=s(krP)(FrP/srP)k2}
  F_{kr}(P) &= F_k(rP)F_r(P)^{k^2}
    &\quad&\text{from Lemma~\ref{lemma:divpolychainrule},} \notag\\
  &= \frac{\s(krP)}{\s(rP)^{k^2}}F_r(P)^{k^2}
    &&\text{from Theorem~\ref{theorem:MTsigma},} \notag\\
  &= \s(krP)\left(\frac{F_r(P)}{\s(rP)}\right)^{k^2}.
\end{alignat}

We claim that the second factor is a $\gp$-adic unit. To see this, we
observe that the composition $\s\circ[r]$ is well-defined on the set
\[
  U_r = \Ecal[r]+\Ecal^f(\Rbar),
\]
which is a $\gp$-adic analytic neighborhood of the $r$-torsion
sections of the scheme~$\Ecal$. Further, since~$\s$ itself has
divisor~$(\Ocal)$ in~$\Ecal^f(\Rbar)$, we see that the divisor of~$\s\circ[r]$
on the set~$U_r$ is given by
\[
  (\s\circ[r])\big|_{U_r} = [r]^*(\Ocal).
\]
On the other hand, the function~$F_r$ has divisor
\[
  (F_r) = [r]^*(\Ocal) - r^2(\Ocal).
\]
Thus
\begin{equation}
  \label{equation:divFronUr}
  \left.\left(\frac{F_r}{\s\circ[r]}\right)\right|_{U_r} = -r^2(\Ocal).
\end{equation}
We have assumed that $r\ge2$, which is equivalent to the assumption
that~$P\notin\Ecal^f$, so~$P$ and~$\Ocal$ do not intersect on the special
fiber of~$\Ecal$. (More formally,~$P$ and~$\Ocal$ determine sections
$s_P,s_\Ocal:\Spec(R)\to\Ecal$, and our assumption ensure that the divisors
$s_P(\Spec(R))$ and $s_\Ocal(\Spec(R))$ do not intersect on the special
fiber~$\Ecal\times_R(R/\gp)$.) It follows from~\eqref{equation:divFronUr}
that $(F_r/\s\circ[r])(P)$ is a $\gp$-adic unit.

Using this fact in~\eqref{equation:Fkr(P)=s(krP)(FrP/srP)k2}, we have
proven that there is a unit $\a\in R^*$ so that
\begin{equation}
  \label{equation:FkrP=skrptimesak2}
  F_{kr}(P) = \s(krP)\cdot \a^{k^2}
  \qquad\text{for all $k\ge1$.}
\end{equation}
The point~$rP$ is in the formal group, so we need to determine the
periodicity properties of~$\s\circ[k]$ on the formal
group~$\Ecal^f(R)$. We use the following well-known elementary result,
whose proof we briefly sketch.

\begin{lemma}
\label{lemma:formalgroupannihilator}
Let $K/\QQ_p$ be a finite extension with ring of integers~$R$, maximal
ideal~$\gp$ and ramification index~$e=\ord_\gp(p)$. Let~$\Gcal/R$ be
a one-parameter formal group.
\begin{parts}
\Part{(a)}
For every $\l\ge1$ there are power series $A_{\l,i}(z)\in zR[[z]]$
for $0\le i\le \l$ so that
\[
  [p^\l]_\Gcal(z) = \sum_{i=0}^\l p^{\l-i}A_{i,\l}\bigl(z^{p^i}\bigr).
\]
\Part{(b)}
Fix $\mu\ge1$, and let~$\l$ be the smallest positive integer with the
property that
\begin{equation}
  \label{equation:deflambda2}
  \min_{0\le i\le \l} \bigl\{ (\l-i)e+p^i \bigr\} \ge \mu.
\end{equation}
Then
\[
  [p^\l]_\Gcal(z) \equiv 0 \pmod{\gp^\mu}
  \qquad\text{for all $z\in\gp$.}
\]
Equivalently, we have
\[
  [p^\l]_\Gcal\bigl(\Gcal(\gp)\bigr) \subset \Gcal(\gp^\mu).
\]
\end{parts}
\end{lemma}
\begin{proof}[Proof of Lemma~\ref{lemma:formalgroupannihilator}]
(a) The multiplication-by-$p$ map on any formal group has the form
\[
  [p](z) = pF(z) + G(z^p)
  \qquad\text{for some $F,G\in R[[z]]$.}
\]
This is most easily proven using the invariant differential, see for
example~\cite[Corollary~IV.4.4]{AEC}. This proves~(a) for $\l=1$.
The general case is then easily proven by induction, using
the formula $[p^{\l+1}](z) = [p]\bigl([p^\l](z)\bigr)$.
\par
In order to prove~(b), we observe that~(a) implies that
\[
  \ord_\gp\bigl([p^\l](z)\bigr)
  \ge\min_{0\le i\le \l} \bigl\{ (\l-i)e+p^i \bigr\}
  \qquad\text{for all $z\in\gp$.}
\]
Then our choice of~$\l$ to satisfy~\eqref{equation:deflambda2} yields
\[
  \ord_\gp\bigl([p^\l](z)\bigr) \ge \mu
  \qquad\text{for all $z\in\gp$.}
\]
This is just another way of saying the $[p^\l](z)\equiv0\MOD{\gp^\mu}$,
which completes the proof of Lemma~\ref{lemma:formalgroupannihilator}.
\end{proof}

We resume the proof of the theorem and we assume that~$\l$ is chosen
as specified in~\eqref{equation:deflambda}, so
Lemma~\ref{lemma:formalgroupannihilator} tells us that
\[
  [p^\l](Q) \equiv 0 \pmod{\gp^\mu}
  \qquad\text{for all $Q\in\Ecal^f(\Rbar)$.}
\]
In particular, this is true for $Q=rP$, where~$P$ was our original
point whose order modulo~$\gp$ is~$r$. Hence
\[
  [p^\l r](P) = [p^\l](rP) \equiv 0 \pmod{\gp^\mu},
\]
so for any $k,j\ge1$ we have
\begin{equation}
  \label{equation:pointperiodicitymodpmu}
  [(k+mp^\l)r](P) = [kr](P) + [m]\bigl([p^\l r](P)\bigr)
  \equiv [kr](P) \pmod{gp^\mu}.
\end{equation}
\par
Substituting this into~\eqref{equation:FkrP=skrptimesak2}, we
find that for all $k,m\ge1$, 
\begin{alignat*}{2}
  F_{(k+mp^\l)r}(P)
  &= \s\bigl((k+mp^\l)rP\bigr) \cdot \a^{(k+mp^\l)^2}  
    &\quad&\text{from \eqref{equation:FkrP=skrptimesak2},} \\
  &\equiv \s(krP)\cdot \a^{(k+mp^\l)^2}  \pmod{\gp^\mu}
    &&\text{from \eqref{equation:pointperiodicitymodpmu},} \\
  &= F_{kr}(P)\cdot  \a^{(k+mp^\l)^2-k^2}
    &\quad&\text{from \eqref{equation:FkrP=skrptimesak2},} \\
  &= F_{kr}(P)\cdot  \bigl(\a^{2+mp^{\l}}\bigr)^{mp^\l}.
\end{alignat*}
In particular, if~$mp^\l$ is a multiple of
\[
  \#(R/\gp^\mu) = \Norm\gp^\mu - \Norm\gp^{\mu-1},
\]
then $\b^{mp^\l}\equiv1\MOD{\gp^\mu}$ for all $\gp$-adic units~$\b$.
However, we observe that taking~$i=0$ in our
definition~\eqref{equation:deflambda} of~$\l$, we have
\text{$\l e+1\ge\mu$}, which implies that~$\Norm\gp^{\mu-1}$ automatically
divides~$p^\l$. Thus it suffices to take~$m$ divisible by
\text{$\Norm\gp-1$}.
\par
We have proven that
if $\ell\ge1$ is any integer satisfying
\[
  p^\l(\Norm\gp-1) \bigm| \ell,
\]
where~$\l$ is chosen to satisfy~\eqref{equation:deflambda}, then  
\[
  F_{(k+\ell)r}(P) \equiv F_{kr}(P) \pmod{\gp^\mu}.
\]
Hence the sequence $(F_{kr}(P)\bmod\gp^\mu)$ is periodic and its period
is as specified in the statement of
Theorem~\ref{theorem:Fnperiodmodptothemu}.
\end{proof}

\begin{remark}
\label{remark:periodicitymodppower}
For $K=\QQ_p$, we always have 
\[
  \min_{0\le i\le \l} \bigl\{ (\l-i)+p^i \bigr\} = \l+1.
\]
Thus the condition~\eqref{equation:deflambda} becomes simply
$\l=\mu-1$, so Theorem~\ref{theorem:Fnperiodmodptothemu} tells
us that
\[
  (F_{kr}(P)\bmod p^\mu)_{k\ge1}
  \quad\text{has period dividing $p^{\mu-1}(p-1)$.}
\]
Taking~$\mu=2$, we recover Shipsey's result, albeit in the context
of division polynomials rather than elliptic divisibility sequences,
and only in the case of good ordinary reduction. 
\par
Continuing with the case $K=\QQ_p$, we consider anew the
formula~\eqref{equation:FkrP=skrptimesak2}, which now says that there
is an~$\a\in\ZZ_p^*$ so that
\begin{equation}
  \label{equation:shprmk1}
  F_{kr}(P) = \s(krP)\cdot \a^{k^2}
  \qquad\text{for all $k\ge1$.}
\end{equation}
We recall that~$rP\in\Ecal^f(R)$, so $z(rP)\equiv0\MOD{p}$.  Thus
\[
  z(krP) = [k](z(rP)) \equiv k\cdot z(rP) \pmod{z(rP)^2},
\]
so in particular, $z(krP)\equiv k\cdot z(rP)\MOD{p^2}$.
Hence with our normalization of the Mazur-Tate $\s$-function, 
it follows that
\begin{equation}
  \label{equation:shprmk2}
  \s(krP) \equiv k\cdot z(rP) \pmod{p^2}.
\end{equation}
Applying~\eqref{equation:shprmk1} and~\eqref{equation:shprmk2} twice,
once with arbitrary~$k$ and once with~$k=1$, we deduce the
simple formula
\begin{equation}
  \label{equation:FkrP=kak2-1FrP}
  F_{kr}(P)   \equiv k\cdot\a^{k^2-1}\cdot F_r(P) \pmod{p^2}.
\end{equation}
This may be compared with Shipsey's formula~\cite[Theorem~3.5.4]{SHIP}
for an elliptic divisibility sequence~$(W_n)$ modulo~$p^2$, for
which she proves (under suitable hypotheses) 
\[
  W_{kr} \equiv k\cdot \b^{k^2-1} \cdot W_r \pmod{p^2}.
\]
(We have simplified Shipsey's formula by observing that
$(-1)^{k+1}$ is equal to $(-1)^{k^2-1}$, so our~$\b$ is the negative 
of Shipsey's~$b$.)
\par
It is interesting to note that an analogous formula for
\text{$F_{kr}(P)\bmod{p^3}$} would necessarily be more complicated,
since it becomes necessary to consider more than the first term of the
power series for~$[k](z)$ and~$\s(z)$.  On the other hand, using the
fact that
\[
  F_{p^{\mu-1}r}(P) \equiv 0 \pmod{p^\mu},
\]
it is possible to give a simple formula for the sequence
\[
  \bigl( F_{kp^{\mu-1}r}(P)  \bmod p^{\mu+1} \bigr)_{k\ge1}
\]
that generalizes~\eqref{equation:FkrP=kak2-1FrP}.
We leave the details to the interested reader.
\end{remark}

\section{Elliptic divisibility sequences}
\label{section:preliminaries}

We are going to use Theorem~\ref{theorem:divpolyconverge} to partially
prove a conjecture concerning the $p$-adic behavior of classical
elliptic divisibility sequences.  Our inspiration for this result, and
indeed the original motivation for much of the work in this paper, is
aptly summarized by the following quote from Morgan Ward's
monograph~\cite[page~33]{Ward1}.

\begin{quotation}
If the least positive residues modulo~$m$ of the successive
values~$U_0,U_1,U_2,\ldots$ of any Lucas function (i.e.,
$U_n=(\a^n-\b^n)/(\a-\b)$) are calculated, the pattern of residues 
exhibits interesting symmetries. These symmetries extend to elliptic
sequences, and find their ultimate explanation in the periodicity of the
second kind of the Weierstrass sigma function.
\end{quotation}

We recall the general definition of a \textit{divisibility
sequence} as being a sequence of integers~$(D_n)_{n\ge0}$ satisfying
\[
  m|n \Longrightarrow D_m|D_n.
\]
A standard example of a divisibility sequence is one of the form
\text{$a^n-1$}, and more generally, divisibility sequences may appear
as linear recurrence sequences such as the Fibonacci sequence.  A
complete classification of divisibility sequences associated to linear
recurrences is given in~\cite{BPvdP}.
\par
It is less clear that there are interesting divisibility sequences
satisfying nonlinear recurrences. The most famous examples of such
sequences are associated to the recursion formula for division
polynomials on  elliptic curves.

\begin{definition}
An \emph{elliptic divisibility sequence} (abbreviated EDS)
is a divisibility sequence  $\Wcal=(W_n)_{n\ge0}$ satisfying the formula
\begin{multline}
  \label{equation:EDSrecursion}
  \W_{m+n}  \W_{m-n} = \W_{m+1}  \W_{m-1}  W_n^2
     - \W_{n+1}  \W_{n-1}  W_m^2 \\
  \text{for all $m\ge n\ge 1$.}
\end{multline}
\end{definition}

The arithmetic properties of elliptic divisibility sequences were
first studied in detail by Morgan Ward~\cite{Ward1,Ward2} in the
1940's, and recently there has been a resurgence of interest in their
study~\cite{CC,Durst,EGW,EPSW,EW1,EW2,SHIP,SW,SS}. (See
also~\cite{Gale,Propp,Rob} for work on the related Somos sequences.)
Ward calls an EDS \emph{proper} if
\[
  W_0=0,\qquad W_1=1,\qquad\text{and}\qquad W_2W_3\ne 0,
\]
and he proves that a proper EDS is associated to a pair~$(E_\Wcal,P_\Wcal)$
consisting of a (possibly singular) elliptic curve and a rational point
$P_\Wcal\in E_\Wcal(\QQ)$. Further, the sequence~$\Wcal$ satisfies a linear
recurrence if and only if the curve~$E_\Wcal$ is singular and is bounded
if and only if~$P_\Wcal$ is a torsion point. 

\begin{remark}
The nonsingularity of the curve~$E_\Wcal$ associated to a proper
elliptic divisibility sequence~$\Wcal$ is equivalent to the nonvanishing
of the \emph{discriminant}
\begin{align}
  \Disc(\Wcal) = {\W_4}W_2^{15} &- {W_3^{3}}W_2^{12} +
  {{3}W_4^{2}}W_2^{10}  - {{{20}\W_4}W_3^{3}}W_2^{7} \notag \\
  &+ {{3}W_4^{3}}W_2^{5} 
    + {{16}W_3^{6}}W_2^{4} 
  + {{{8}W_4^{2}}W_3^{3}}W_2^{2} +  {W_4^{4}}.
  \label{equation:Disc(W)}
\end{align}
This is essentially the discriminant of the curve~$E_\Wcal$
(cf{.} Ward~\cite[equation (19.3)]{Ward1}).
See~\cite{Ward1} or~\cite[Appendix]{SS} for additional formulas
describing~$E_\Wcal$ and~$P_\Wcal$.
\end{remark}

\begin{remark}
Except in some degenerate cases, the definition of an~EDS forces
$W_0=0$ (put $m=n$ in~\eqref{equation:EDSrecursion}) and $W_1=\pm1$
(put $n=1$).  As noted above, a proper EDS with $\Disc(\Wcal)=0$
satisfies a linear recurrence.  See~\cite{SHIP,Ward1} for details and
a complete description of nonproper EDS.
\end{remark}

\begin{definition}
We call~$\Wcal$ a
\emph{general elliptic divisibility sequences} if it satisfies:
\begin{enumerate}
\item
$\Wcal$ is proper;
\item
$E_\Wcal$ is a nonsingular elliptic curve, or equivalently,
$\Disc(\Wcal)\ne0$;
\item
$P_\Wcal$ is a point of infinite order in $E_\Wcal(\QQ)$, or
equivalently,~$\Wcal$ is unbounded.
\end{enumerate}
From our earlier remarks, these are the most interesting EDS.
(Note that this definition differs somewhat from Ward's terminology.)
\end{definition}

\begin{example}
\label{example:EDS37}
The simplest general elliptic divisibility sequence is
the sequence
\begin{multline}
  1, 1, -1, 1, 2, -1, -3, -5, 7, -4, -23, 29, 59, 129, -314, \\
  -65, 1529, -3689, -8209, -16264,\dots
\end{multline}
It is associated to the generator~$P=(0,0)$ of the Mordell-Weil group
on the elliptic curve \text{$y^2+y=x^3-x$}  of conductor~37.
\end{example}

An elliptic divisibility sequence~$(W_n)$ is required to satisfy 
the recursion~\eqref{equation:EDSrecursion} for all $m\ge n\ge1$.
It is easy to check that it suffices that~$(W_n)$ satisfy
the two relations
\begin{align}
  \label{equation:EDSrecursionodd}
  \W_{2n+1} &=\W_{n+2}W_n^3 - \W_{n-1}W_{n+1}^3, \\
  \label{equation:EDSrecursioneven}
  \W_{2n}\W_2 &= \W_n\left(\W_{n+2}W_{n-1}^2 -\W_{n-2}W_{n+1}^2\right).
\end{align}
In particular, a proper EDS is determined by the values
of~$W_2,W_3,W_4$. Further, a triple~$W_2,W_3,W_4$ with $W_2W_3\ne0$
gives an EDS if and only if~$W_2|W_4$. (See~\cite{Ward1}.)  
Additional material on elliptic divisibility sequences may be found in
\cite{Durst,EGW,EPSW,EW1,EW2,SHIP,SW,SS,Ward1,Ward2}.

\section{Elliptic divisibility sequences and elliptic functions}

We recall Ward's fundamental result relating elliptic divisibility
sequences to values of elliptic functions.

\begin{theorem}
\label{theorem:EDS=sigma}
Let $(W_n)$ be a general elliptic divisibility sequence. Then
there is a lattice $L\subset\CC$ and a complex number $\xi\in\CC$
such that
\[
  W_n = \psi_n(\xi,L) = \frac{\sigma(n\xi,L)}{\sigma(\xi,L)^{n^2}}
  \qquad\text{for all $n\ge1$,}
\]
where~$\psi_n(\xi,L)$ and~$\sigma(\xi,L)$ are, respectively, the
{\upshape(}analytic{\upshape)} $n$-division polynomial and the
Weierstrass $\s$-function associated to the lattice~$L$.
\par
Further, the modular invariants~$g_2(L)$ and~$g_3(L)$ associated to
the lattice~$L$ and the Weierstrass values~$\wp(\xi,L)$
and~$\wp'(\xi,L)$ associated to the point~$\xi$ on the elliptic
curve~$\CC/L$ are in the field~$\QQ$.  {\upshape(}More
precisely,~$g_2(L),g_3(L),\wp(\xi,L),\wp'(\xi,L)$ are given as
rational expressions in~$\QQ(W_2,W_3,W_4)$.{\upshape)}
\end{theorem}
\begin{proof}
See Ward~\cite[Theorems 12.1 and~19.1]{Ward1}. The rational
expressions for~$g_2$ and~$g_3$ in~$\QQ(W_2,W_3,W_4)$ are given by
\cite[equations~13.6 and~13.7]{Ward1}, and the rational expressions
for~$\wp(\xi,L)$ and~$\wp'(\xi,L)$ are given by \cite[equations~13.5
and~13.1]{Ward1}.  See also~\cite[Appendix]{SS} for simplified
formulas.
\end{proof}

We also mention that Ward proves a partial converse to
Theorem~\ref{theorem:EDS=sigma}.

\begin{theorem}
Let~$L$ be a lattice with $g_2(L),g_3(L)\in\QQ$ and let~$\xi\in\CC$
satisfy $\wp(\xi,L),\wp'(\xi,L)\in\QQ$. Then there is a
constant~$c\in\QQ^*$ so that the sequence $(c^{n^2-1}\psi_n(\xi,L))$
is an elliptic divisibility sequence.
\end{theorem}
\begin{proof}
This is proven in \cite[Theorem~21.4]{Ward1}.  See also~\cite{SHIP,SW}.
\end{proof}

We reformulate Ward's result so that it is entirely in terms of
rational numbers.

\begin{proposition}
\label{proposition:EDS=Fn}
Let $\Wcal=(W_n)$ be a nonsingular elliptic divisibility sequence, and
let~$E_\Wcal/\QQ$ and~$P_\Wcal\in E(\QQ)$ be the associated elliptic
curve and rational point.  Fix a minimal
Weierstrass equation for~$E_\Wcal$, and let~$F_n$ be the normalized
$n$-division polynomial on~$E_\Wcal$
{\upshape(}Definition~\ref{definition:ndivpoly}{\upshape)}.  Then
there is a constant $\g\in\QQ^*$ so that
\[
  W_n = \g^{n^2-1} F_n(P_\Wcal)\qquad\text{for all $n\ge1$.}
\]
Further, the denominator of~$\g$ is divisible only by primes of bad
reduction of~$P_\Wcal$, i.e., primes~$p$ of bad reduction
for~$E_\Wcal$ at which~$P_\Wcal\bmod p$ is the singular point
on~$E_\Wcal\bmod p$.
\end{proposition}
\begin{proof}
We use Theorem~\ref{theorem:EDS=sigma} to choose a lattice~$L$ and
complex number~$\xi\in\CC$ so that the given EDS has the
form $W_n = \psi_n(\xi,L)$ for all $n\ge1$. Let
$\F:\CC/L\to E(\CC)$ be an isomorphism. Then Lemma~\ref{lemma:Fvspsi}
tells us that there is a constant~$\g\in\CC^*$ so that
\[
  F_n(\F(\z)) = \g^{1-n^2}\psi_n(\z,L)
  \qquad\text{for all $\z\in\CC$ and all $n\ge1$.}
\]
Substituting~$\z=\xi$, so~$P_\Wcal=\F(\xi)$, we obtain
\[
  F_n(P_\Wcal) = \g^{1-n^2} \psi_n(\xi,L) = \g^{1-n^2} W_n
  \qquad\text{for all $n\ge1$.}
\]
Putting~$n=2$ and~$n=3$, we find that~$\g^3$ and~$\g^8$ are in~$\QQ$,
and hence that~$\g\in\QQ$, which completes the proof that
$W_n=\g^{n^2-1}F_n(P_\Wcal)$ with $\g\in\QQ^*$.  
\par
Since~$W_n\in\ZZ$, we see that
\begin{align}
  \label{equation:ordpFnP}
 \ord_p(F_n(P_\Wcal)) &= \ord_p{W_n} - (n^2-1)\ord_p(\g) \notag\\
 &\ge  -(n^2-1)\ord_p(\g)   \qquad\text{for all $n\ge1$. }
\end{align}
Let~$p$ be a prime such that \text{$P_\Wcal\bmod p$} is nonsingular.
Then~\text{$P_\Wcal\bmod p$} cannot have order both~2 and~3, so
at least one of~$F_2(P_\Wcal)$ and~$F_3(P_\Wcal)$ is nonzero modulo~$p$.
Then~\eqref{equation:ordpFnP} says that either
\[
  0 = \ord_p(F_2(P_\Wcal)) \ge -3\ord_p(\g)
  \text{ or }
  0 = \ord_p(F_3(P_\Wcal)) \ge -3\ord_p(\g),
\]
and hence $\ord_p(\g)\ge0$.
\end{proof}

\section{A $p$-adic limit of elliptic divisibility sequences}
\label{section:EDSpadiclimit}

In this section we apply our results on division polynomials
to partially prove the following conjecture about elliptic divisibility
sequences.

\begin{conjecture}
\label{conjecture:EDSpadicconverge}
Let $\Wcal=(W_n)_{n\ge0}$ be an elliptic divisibility sequence
and let~$p$ be a prime. Then there is an exponent $N=N_p\ge1$
so that for every $m\ge1$, the limit
\[
  \lim_{k\to\infty} W_{mp^{kN}}
\]
converges in~$\ZZ_p$ to a number that is algebraic over~$\QQ$.
\end{conjecture}

We are able to prove this conjecture for ``most'' elliptic divisibility
sequences and for ``most'' primes.

\begin{theorem}
\label{theorem:EDSpadicconverge}
Let $\Wcal=(W_n)_{n\ge0}$ be a general elliptic divisibility sequence.
Let~$E_\Wcal$ be the associated elliptic curve, given by a minimal
Weierstrass equation over~$\QQ$, let~$P_\Wcal\in E_\Wcal(\QQ)$ be the
associated rational point, and let~$S_\Wcal$ be the set of primes~$p$
satisfying any one of the following conditions:
\par\noindent
\vbox{%
\begin{itemize}
\item
$p=2$.
\item
$P_\Wcal\equiv\Ocal\MOD{p}$.
\item
$P_\Wcal\bmod p$ is a singular point on $E_\Wcal\bmod p$.
\item
$E_\Wcal\bmod p$ is supersingular.
\end{itemize}
}
\par\noindent
Then Conjecture~\ref{conjecture:EDSpadicconverge} is true for~$\Wcal$
for all primes $p\notin S_\Wcal$, that is, 
there is an exponent $N=N_p\ge1$
so that for every $m\ge1$, the limit
\[
  \lim_{k\to\infty} W_{mp^{kN}}
\]
converges in~$\ZZ_p$ to a number that is algebraic over~$\QQ$.
\par
In particular, if~$E_\Wcal$ does not have complex multiplication, then
Conjecture~\ref{conjecture:EDSpadicconverge} is true for almost all
primes in the sense of density, since~$S_\Wcal$ has density~$0$.
{\upshape(}If~$E_\Wcal$ has~CM, then~$S_\Wcal$ has
density~$\frac12$.{\upshape)}
\end{theorem}

\begin{proof}
We use Proposition~\ref{proposition:EDS=Fn} to find a~$\g\in\QQ^*$ so
that
\begin{align} 
  \label{equation:Wn=gn2-1FnP}
  W_n = \g^{n^2-1}F_n(P_\Wcal)
  \qquad\text{for all $n\ge1$.}
\end{align}
Proposition~\ref{proposition:EDS=Fn} also tells us
that~$\ord_p(\g)\ge0$, since we have assumed that \text{$P_\Wcal\bmod
p$} is nonsingular. On the other hand, our assumption
that \text{$P_\Wcal\not\equiv\Ocal\MOD{p}$} implies
that $\ord_p(F_n(P_\Wcal))\ge0$, so if $\ord_p(\g)>0$, then
\[
  \ord_p(W_n) \ge (n^2-1)\ord_p(\g) \longrightarrow \infty
  \quad\text{as $n\to\infty$.}
\]
Thus if $\ord_p(\g)>0$, then $\lim_{n\to\infty} W_n=0$ in~$\ZZ_p$.
\par
We are thus reduced to the case that~$\g\in\ZZ_p^*$. Then for 
any $m\ge1$,
\[
  \lim_{k\to\infty} \g^{mp^k-1}
  = \chi(\g)^m\g^{-1}\quad\text{in $\ZZ_p$,}
\]
where~$\chi(\g)\in\bfmu_{p-1}$ is the value of the Teichm\"uller 
character. Using~\eqref{equation:Wn=gn2-1FnP}, it thus suffices
to prove there is a power~$q=p^N$ so that for every $m\ge1$,
\[
  \lim_{k\to\infty} F_{mq^k}(P_\Wcal)
  \quad\text{converges in $\ZZ_p$}
\]
and is algebraic over~$\QQ$. This follows immediately from
Theorem~\ref{theorem:divpolyconverge}, since our choice of
the set~$S_\Wcal$ was designed to ensure that
Theorem~\ref{theorem:divpolyconverge} is applicable to
every prime~$p\notin S_\Wcal$.
\par
This completes the proof of Theorem~\ref{theorem:EDSpadicconverge}
except for the final statements about densities. We note that the
first three conditions specifying primes in~$S_\Wcal$ only include
finitely many primes. The fourth condition, that~$E_\Wcal\MOD{p}$ be
supersingular, is more serious.  However, Serre~\cite{SE1} has proven
that for any fixed non-CM elliptic curve~$E/\QQ$, almost all primes
give ordinary reduction.  (More precisely, the number of supersingular
primes less than~$X$ is~$O(X^{3/4})$, see~\cite{EL2,SE2}.) On the
other hand, Elkies has shown that there are infinitely many primes of
supersingular reduction~\cite{EL1}, so unfortunately our set~$S_\Wcal$
is always infinite.  Finally, it is well known that an elliptic curve
with~CM has supersingular reduction at precisely the primes that are
inert in its CM~field.
\end{proof}

\begin{remark}
Theorem~\ref{theorem:divpolyconverge} tells us when the limit in
Theorem~\ref{theorem:EDSpadicconverge} is nonzero.  Let~$r_p$ be the
order of~$P_\Wcal$ in~$E_\Wcal(\FF_p)$, and let $q=p^N$ be the power
of~$p$ in Theorem~\ref{theorem:EDSpadicconverge}.  Note that Hasse's
theorem~\cite[V.1.1]{AEC} implies that $r_p\le(\sqrt p+1)^2$, so in
particular, $r_p<p^2$ under our assumption that $p\ne2$.
\par
Writing $W_n=\g^{n^2-1}F_n(P_\Wcal)$ as in~\eqref{equation:Wn=gn2-1FnP},
we see that if~$p|\g$, then $p|W_n$ for all~$n\ge2$, in which case
$W_n\to0$ in~$\ZZ_p$. On the other hand, if~$p\notdivide\g$, 
then Theorem~\ref{theorem:divpolyconverge} implies that
the limit in Theorem~\ref{theorem:EDSpadicconverge} is zero if
and only if the $p$-free part of~$r_p$ divides~$m$ Hence
\[
  \lim_{k\to\infty} W_{mp^{kN}} = 0
  \quad\text{if and only if either}\quad
  \begin{cases}
    \text{$p|W_n$ for all $n\ge2$, or} \\
    \text{$r_p|mp$.} \\
  \end{cases}
\]
\end{remark}


\begin{acknowledgement}
The author would like to thank Noam Elkies, Graham Everest, Barry
Mazur, Rachel Shipsey, Nelson Stephens, and Thomas Ward for helpful
correspondence during the preparation of this paper.
\end{acknowledgement}




\begin{thebibliography}{99}

\itemsep=\smallskipamount

\bibitem{BPvdP}
J.P. B\'ezivin, A. Peth\"o, A.J. van der Poorten, A full
characterization of divisibility sequences, \emph{Amer. J. of
Math.} \textbf{112} (1990), 985--1001.


\bibitem{CH}
J. Cheon, S. Hahn, Explicit valuations of division polynomials of an
elliptic curve, \emph{Manuscripta Math.} \textbf{97} (1998), 319--328.

\bibitem{CC}
D.V. Chudnovsky, G.V. Chudnovsky, Sequences of numbers generated
by addition in formal groups and new primality and factorization tests,
\emph{Advances in Applied Mathematics} \textbf{7} (1986), 385--434.

\bibitem{Durst}
L.K. Durst, The apparition problem for equianharmonic divisibility,
Proc. Nat. Acad. Sci. U. S. A. \textbf{38} (1952), 330--333.

\bibitem{EGW}
M. Einsiedler, G. Everest, T. Ward, Primes in elliptic divisibility
sequences, \emph{LMS J. Comput. Math.} \textbf{4} (2001), 1--13,
{electronic}.

\bibitem{EL1}
N. Elkies, 
{The existence of infinitely many supersingular primes for
every elliptic curve over {${\QQ}$}},
\emph {Invent. Math.} \textbf{89}, (1987), 561--567.

\bibitem{EL2}
\bysame
{Distribution of supersingular primes},
{Journ\'ees Arithm\'etiques, 1989 (Luminy, 1989)},
\emph{Ast\'erisque} \textbf{198-200} (1991),
{127--132 (1992)}.

\bibitem{EPSW}
G. Everest, A. van der Poorten, I. Shparlinski, T. Ward,
\emph{Recurrence sequences}, Mathematical Surveys and Monographs
104, AMS, Providence, RI, 2003.


\bibitem{EW1}
G. Everest, T. Ward, Primes in divisibility sequences,
\emph{Cubo Mat. Educ.} \textbf{3} (2001), 245--259.

\bibitem{EW2}
\bysame
The canonical height of an algebraic point on an elliptic curve,
\emph{New York J. Math.} 6 (2000), 331-342, (electronic).

\bibitem{Gale}
D. Gale, 
The Strange and Surprising Saga of the Somos Sequences,
\emph{Mathematical Intelligencer} \textbf{13}, (1991) 40--42, 49--50




\bibitem{MT}
B. Mazur, J. Tate, The $p$-adic sigma function, \emph{Duke Math. J.}
\textbf{62} (1991), 663--688.


\bibitem{Propp}
J. Propp, 
The Somos Sequence Site, \\
\texttt{http://www.math.wisc.edu/\~{}propp/somos.html}.

\bibitem{Rob}
R. Robinson, Periodicity of Somos
sequences, \emph{Proc. Amer. Math. Soc.} \textbf{116} (1992), 613--619.

\bibitem{Schoof1}
R. Schoof, Elliptic curves over finite fields and the computation of
square roots mod~$p$, \emph{Math. Comp.} \textbf{44} (1985), 483--494.

\bibitem{Schoof2}
\bysame
Counting points on elliptic curves over finite fields, Les
Dix-huiti\`emes Journ\'ees Arithm\'etiques (Bordeaux,
1993), \emph{J. Th\'eor. Nombres Bordeaux} \textbf{7} (1995), 219--254.

\bibitem{SE1}
J.-P. Serre,
\emph{Abelian {$l$}-adic representations and elliptic curves},
{W. A. Benjamin, Inc., New York-Amsterdam}, {1968}.

\bibitem{SE2}
\bysame 
{Quelques applications du th\'eor\`eme de densit\'e de
{C}hebotarev},
\emph{Inst. Hautes \'Etudes Sci. Publ. Math.},
\textbf {54} (1981), {323--401}.

\bibitem{SHIP}
R. Shipsey, 
Elliptic divisibility sequences, Ph.D. thesis, Goldsmith's College
(University of London), 2000.

\bibitem{AEC}
J.H. Silverman, \emph{The arithmetic of elliptic curves}, GTM~106, 
Springer-Ver\-lag, New York, 1986.

\bibitem{ATAEC}
\bysame
\emph{Advanced topics in the arithmetic of
elliptic curves}, GTM~151, Springer-Ver\-lag, New York, 1994.


\bibitem{SS}
J.H. Silverman, N. Stephens, The sign of an elliptic divisibility
sequence, preprint 2004,  (arXiv:mathNT/0402415).

\bibitem{SW}
C.S. Swart, 
Elliptic divisibility sequences, Ph.D. thesis, Royal Holloway
(University of London), 2003.


\bibitem{Ward1}
M. Ward, Memoir on elliptic divisibility sequences, \emph{Amer. J. Math.}
\textbf{70} (1948), 31--74.

\bibitem{Ward2}
M. Ward, The law of repetition of primes in an elliptic divisibility
sequence, \emph{Duke Math. J.} \textbf{15} (1948), 941--946.

\bibitem{WW}
E.T. Whittaker, G.N. Watson, \emph{A course in modern analysis}, 
Cambridge Univ. Press, Cambridge, 4th ed., 1927.

\end{thebibliography}
\end{document}